\documentclass{article}
\usepackage{amsfonts,amsmath,amsxtra}
\usepackage{latexsym}
\usepackage{amssymb}

\usepackage{comment}

\def\hybrid{\topmargin 0pt      \oddsidemargin 0pt
        \headheight 0pt \headsep 0pt
        \textwidth 16.5cm
        \textheight 23cm
        \marginparwidth 0.0in
        \parskip 5pt plus 1pt   \jot = 1.5ex}
\catcode`\@=11
\def\marginnote#1{}
\newcount\hour
\newcount\minute
\newtoks\amorpm
\hour=\time\divide\hour by60 \minute=\time{\multiply\hour by60
\global\advance\minute by-\hour}
\edef\standardtime{{\ifnum\hour<12 \global\amorpm={am}%
        \else\global\amorpm={pm}\advance\hour by-12 \fi
        \ifnum\hour=0 \hour=12 \fi
      \number\hour:\ifnum\minute<10 0\fi\number\minute\the\amorpm}}
\edef\militarytime{\number\hour:\ifnum\minute<10 0\fi\number\minute}

\def\draftlabel#1{{\@bsphack\if@filesw {\let\thepage\relax
   \xdef\@gtempa{\write\@auxout{\string
      \newlabel{#1}{{\@currentlabel}{\thepage}}}}}\@gtempa
   \if@nobreak \ifvmode\nobreak\fi\fi\fi\@esphack}
        \gdef\@eqnlabel{#1}}
\def\@eqnlabel{}
\def\@vacuum{}
\def\draftmarginnote#1{\marginpar{\raggedright\scriptsize\tt#1}}

\def\draft{\oddsidemargin -0.1truein
        \def\@oddfoot{\sl preliminary draft \hfil
        \rm\thepage\hfil\sl\today\quad\militarytime}
        \let\@evenfoot\@oddfoot \overfullrule 3pt
        \let\label=\draftlabel
        \let\marginnote=\draftmarginnote
\def\@eqnnum{{\rm (\theequation)}
\rlap{\kern\marginparsep\tt\@eqnlabel}%
\global\let\@eqnlabel\@vacuum}  }

\newfont{\Bbbb}{msbm7 scaled 1\@ptsize00}
\newcommand{\zs}{\raise-1pt\hbox{$\mbox{\Bbbb Z}$}}

scaled 1\@ptsize00

\font\sevenmsa=msam6 
\newfam\msafam
\textfont\msafam=\sevenmsa
\def\hexnumber@#1{\ifnum#1<10 \number#1\else
\ifnum#1=10 A\else\ifnum#1=11 B\else\ifnum#1=12 C\else \ifnum#1=13
D\else\ifnum#1=14 E\else\ifnum#1=15 F\fi\fi\fi\fi\fi\fi\fi}
\def\msa@{\hexnumber@\msafam}
\def\llcorner{\delimiter"4\msa@78\msa@78 }
\def\lrcorner{\delimiter"5\msa@79\msa@79 }
\mathchardef\blacktriangleright="3\msa@49
\mathchardef\blacktriangleleft="3\msa@4A \font\tenmsb=msbm10 scaled
1\@ptsize00
\newfam\msbfam
\textfont\msbfam=\tenmsb \scriptfont\msbfam=\tenmsb

\newdimen\Squaresize \Squaresize=14pt
\newdimen\Thickness \Thickness=0.5pt

\def\Square#1{\hbox{\vrule width \Thickness
   \vbox to \Squaresize{\hrule height \Thickness\vss
      \hbox to \Squaresize{\hss#1\hss}
   \vss\hrule height\Thickness}
\unskip\vrule width \Thickness} \kern-\Thickness}

\def\Vsquare#1{\vbox{\Square{$#1$}}\kern-\Thickness}

\def\numberbysection{\@addtoreset{equation}{section}
        \def\theequation{\thesection.\arabic{equation}}}
\numberbysection

\renewcommand{\theequation}{\thesection.\arabic{equation}}
\def\titlepage{\@restonecolfalse\if@twocolumn\@restonecoltrue\onecolumn
     \else \newpage \fi \thispagestyle{empty}\c@page\z@
        \def\thefootnote{\fnsymbol{footnote}} }

\def\endtitlepage{\if@restonecol\twocolumn \else  \fi
        \def\thefootnote{\arabic{footnote}}
        \setcounter{footnote}{0}}  
\relax

\hybrid
\makeatletter
\newdimen\normalarrayskip            
\newdimen\minarrayskip               
\normalarrayskip\baselineskip \minarrayskip\jot
\newif\ifold             \oldtrue            \def\new{\oldfalse}
\def\arraymode{\ifold\relax\else\displaystyle\fi}
\def\eqnumphantom{\phantom{(\theequation)}} 
\def\@arrayskip{\ifold\baselineskip\z@\lineskip\z@
     \else
     \baselineskip\minarrayskip\lineskip1\baselineskip\fi}

\def\@arrayclassz{\ifcase \@lastchclass \@acolampacol \or
\@ampacol \or \or \or \@addamp \or
   \@acolampacol \or \@firstampfalse \@acol \fi
\edef\@preamble{\@preamble
  \ifcase \@chnum
     \hfil$\relax\arraymode\@sharp$\hfil
     \or $\relax\arraymode\@sharp$\hfil
     \or \hfil$\relax\arraymode\@sharp$\fi}}

\def\@array[#1]#2{\setbox\@arstrutbox=\hbox{\vrule
     height\arraystretch \ht\strutbox
     depth\arraystretch \dp\strutbox
width\z@}\@mkpream{#2}\edef\@preamble{\halign \noexpand\@halignto
\bgroup \tabskip\z@ \@arstrut \@preamble \tabskip\z@ \cr}%
\let\@startpbox\@@startpbox \let\@endpbox\@@endpbox
  \if #1t\vtop \else \if#1b\vbox \else \vcenter \fi\fi
  \bgroup \let\par\relax
  \let\@sharp##\let\protect\relax
  \@arrayskip\@preamble}
\def\eqnarray{\stepcounter{equation}%
              \let\@currentlabel=\theequation
              \global\@eqnswtrue
              \global\@eqcnt\z@
              \tabskip\@centering              
              \let\\=\@eqncr
              $$%
            \halign to \displaywidth  \bgroup
             \eqnumphantom \@eqnsel
      \hskip\@centering                               
    $\displaystyle  \tabskip\z@ {##}$%
    &\global\@eqcnt\@ne \hskip 2\arraycolsep
         $ \displaystyle  \arraymode{##}$\hfil
    &\global\@eqcnt\tw@ \hskip 2\arraycolsep
         $\displaystyle\tabskip\z@{##}$\hfil
         \tabskip\@centering
    &{##}\tabskip\z@\cr}
\makeatother

\newtheorem{te}{Theorem}[section]
\newtheorem{de}{Definition}[section]
\newtheorem{prop}{Proposition}[section]           
\newtheorem{cor}{Corollary}[section]
\newtheorem{lem}{Lemma}[section]

\newcommand\bqa{\begin{eqnarray}}
\newcommand\eqa{\end{eqnarray}}
\def\be{\begin{eqnarray}\new\begin{array}{cc}}
\def\ee{\end{array}\end{eqnarray}}

\def\beq{\begin{equation}}
\def\eeq{\end{equation}}
\def\bse{\begin{subequations}}                
\def\ese{\end{subequations}}
\def\bp{\begin{pmatrix}}
\def\ep{\end{pmatrix}}

\def\proof{\noindent {\it Proof}. }

\def\stack#1#2{\raise0.7pt\hbox{$\mathrel{\mathop{#2}\limits^{#1}}$}}
\def\tr{\triangleright}
\def\tl{\triangleleft}
\def\sem{\mathsurround=0pt \raise1pt
\hbox{$\scriptscriptstyle>\!\!$}\:\!\!\tl}
\def\mes{\mathsurround=0pt \tr\!\:\!\raise0.8pt
\hbox{$\scriptscriptstyle\!\!<$}\,}
\def\]{\mathsurround=0pt ]\raise-2pt\hbox{$_\ast$}}

\def\<{\langle}
\def\>{\rangle}

\def\we{\raise-1pt\hbox{$\,\stackrel{\wedge}{,}\,$}}
\def\tr{{\rm tr}\,}

\newcounter{pac}[section]

\newcounter{pacc}[subsection]

0
\begin{document}

\setcounter{pac}{0}
\setcounter{footnote}0

\begin{center}

\phantom.
\bigskip
{\Large\bf On explicit realization of algebra of complex powers of generators of $U_{q}(\mathfrak{sl}(3))$}

\vspace{1cm}

\bigskip\bigskip

{\large  Pavel Sultanich \footnote {E-mail:  sultanichp@gmail.com}},\\
\bigskip
{\it Moscow Center for Continuous Mathematical Education, 119002, Bolshoy Vlasyevsky Pereulok 11, Moscow, Russia
}\\
\bigskip

\end{center}




\begin{abstract}
\noindent

In this note  we prove an integral identity involving complex powers of generators of quantum group $U_{q}(\mathfrak{sl}(3))$ considered as certain positive operators in the setting of positive principal series representations. This identity represents a continuous analog of one of the Lusztig's relations between divided powers of generators of quantum groups, which play an important role in the study of irreducible modules \cite{Lu 1}. We also give definitions of arbitrary functions of $U_{q}(\mathfrak{sl}(3))$ generators and give another proofs for some of the known results concerning positive principal series representations of $U_{q}(\mathfrak{sl}(3))$.

\end{abstract}
\vspace{5 mm}

\section{Introduction}

The notion of modular double of a quantum group $U_{q}(\mathfrak{g})$ plays an important role in different areas of mathematical physics such as Liouville theory \cite{PT},\cite{FKV},  relativistic Toda model \cite{KLSTS} and others. It was introduced by Faddeev in \cite{F2} who noticed that certain representations of a quantum group $U_{q}(\mathfrak{sl}(2))$, $q = e^{\pi\imath b^{2}}$ have a remarkable duality under $b\leftrightarrow b^{-1}$ and proposed to consider instead of single quantum group enlarged object generated by two sets of generators $K$, $E$, $F\in U_{q}(\mathfrak{sl}(2))$ and $\tilde{K}$, $\tilde{E}$, $\tilde{F}\in U_{\tilde{q}}(\mathfrak{sl}(2))$, $\tilde{q} = e^{\pi\imath b^{-2}}$. In \cite{BT} it was shown that in a special class of representations of modular double the rescaled generators defined by
$K$, $\mathcal{E} = -\imath(q-q^{-1})E$, $\mathcal{F} = -\imath(q-q^{-1})F$ of $U_{q}(\mathfrak{sl}(2))$ are positive operators. This allows one to use functional calculus and consider arbitrary functions of them. Moreover, the generators of dual group $\tilde{K}$, $\tilde{\mathcal{E}} = -\imath(\tilde{q}-\tilde{q}^{-1})\tilde{E}$, $\tilde{\mathcal{F}} = -\imath(\tilde{q}-\tilde{q}^{-1})\tilde{F}$  are expressed as non-integer powers of the original generators
\begin{equation}
\tilde{K} = K^{b^{-2}},
\end{equation}
\begin{equation}
\tilde{\mathcal{E}} = \mathcal{E}^{b^{-2}},
\end{equation}
\begin{equation}
\tilde{\mathcal{F}} = \mathcal{F}^{b^{-2}}.
\end{equation}
These relations were called transcendental relations. This kind of representations, admitting the transcendental relations, has been generalized to higher ranks \cite{FrIp}, \cite{Ip2} and has been called positive principal series representations. Introduction of particular non-integer powers of generators of quantum group naturally leads to consideration of arbitrary powers of generators.  Thus, the modular double becomes a discrete subalgebra in the algebra generated by arbitrary powers of generators $K^{\imath p}$, $\mathcal{E}^{\imath s}$, $\mathcal{F}^{\imath t}$.

In \cite{Lu 1}, eq.(4.1a)-eq.(4.1j) Lusztig summarized the relations between the divided powers of generators of $U_{q}(\mathfrak{g})$ for simply-laced $\mathfrak{g}$. He used these identities in the study of finite-dimensional modules of $U_{q}(\mathfrak{g})$ in the case where $q$ is a root of unity.
So in the study of the algebra of arbitrary complex powers of quantum group generators, the question of the generalization of these relations arises. Some of the relations were found in \cite{Ip1}, eq.(6.16), eq.(6.17). Another integral relation which is a generalization of Kac's identity \cite{Lu 1}, eq.(4.1a) appeared in \cite{Su1} and was proved in the explicit representation for the case of $U_{q}(\mathfrak{sl}(2))$ in \cite{Su2}. To write it down explicitly, let $G_{b}(x)$ be quantum dilogarithm \cite{F1} which is a special function playing an important role in the study of algebra of complex powers of generators of $U_{q}(\mathfrak{g})$. Its properties will be outlined in Section 2. Let $K_{j} = q^{H_{j}}$, $\mathcal{E}_{j}$, $\mathcal{F}_{j}$ be $U_{q}(\mathfrak{g})$ generators which are assumed to be positive operators so that the functions of them are defined.  Let continuous analogs of divided powers be defined by $A^{(\imath s)}_{i} = G_{b}(-\imath bs)A^{\imath s}$.  Explicit expressions for the powers of operators under consideration will be given later.  Then the generalized Kac's identity reads

\begin{equation}\begin{split}
\mathcal{E}_{j}^{(\imath s)}\mathcal{F}_{j}^{(\imath t)} = \int\limits_{\mathcal{C}} d\tau e^{\pi bQ\tau}\mathcal{F}_{j}^{(\imath t+\imath\tau)}K_{j}^{-\imath \tau}
\frac{G_{b}(\imath b\tau)G_{b}(-bH_{j} + \imath b(s+t+\tau))}{G_{b}(-bH_{j}+\imath b(s+t+2\tau))}\mathcal{E}_{j}^{(\imath s + \imath \tau)},
\end{split}\end{equation}
where the contour $\mathcal{C}$ goes slightly above the real axis but passes below the pole at $\tau = 0$.
In this note we prove this identity for the case of positive principal series representations of $U_{q}(\mathfrak{sl}(3))$.

The paper is organized as follows. In Section 2, we recall the definition of quantum group $U_{q}(\mathfrak{g})$ and outline the definition and basic properties of the quantum dilogarithm $G_{b}(x)$ and related function $g_{b}(x)$. In Section 3 we recall the construction of arbitrary functions of generators and generalized Kac's identity in the case of positive principal series representations of $U_{q}(\mathfrak{sl}(2))$.The main result of the paper is formulated in Theorem 4.1 in Section 4. We define arbitrary functions of $U_{q}(\mathfrak{sl}(3))$ generators in the positive principal series representations. We prove the generalized Kac's identity using unitary transform interwining the formulas of functions of generators of $U_{q}(\mathfrak{sl}(2))_{i}$ subalgebra corresponding to simple root $i$, with the formulas for $U_{q}(\mathfrak{sl}(2))$, defined in Section 3. This calculation also represents another proof of Theorem 4.7 in \cite{Ip1} which states that positive principal series representation of $U_{q}(\mathfrak{sl}(3))$ decomposes into direct integral of positive principal series representations of its $U_{q}(\mathfrak{sl}(2))$ subalgebra corresponding to each simple root.

{\bf Acknowledgements:} The research was supported by    RSF  (project 16-11-10075). I am grateful to A.A.Gerasimov and D.R.Lebedev  for helpful discussions and interest in this work.

\section{Preliminaries}

We start with the definition of quantum groups following \cite{ChPr},\cite{Lu book}.
Let $(a_{ij})_{1\le i,j\le r}$ be Cartan matrix of semisimple Lie algebra $\mathfrak{g}$ of rank $r$. Let $\mathfrak{b}_{\pm}\subset \mathfrak{g}$ be opposite Borel subalgebras. For simplicity let us restrict ourselves to the simply-laced case $a_{ii} = 2$, $a_{ij} = a_{ji} = \{0,-1\}$, $i\ne j$. Let $U_{q}(\mathfrak{g})$ $(q = e^{\pi\imath b^{2}}$, $b^{2}\in \mathbb{R}\setminus \mathbb{Q})$ be the quantum group with generators $E_{j}$, $F_{j}$, $K_{j} = q^{H_{j}}$, $1\le j \le r$ and relations
\begin{equation}
    K_{i}K_{j} = K_{j}K_{i},
\end{equation}
\begin{equation}
    K_{i}E_{j} = q^{a_{ij}}E_{j}K_{i},
\end{equation}
\begin{equation}
    K_{i}F_{j} = q^{-a_{ij}}F_{j}K_{i},
\end{equation}
\begin{equation}
    E_{i}F_{j} - F_{j}E_{i} = \delta_{ij}\frac{K_{i} - K_{i}^{-1}}{q-q^{-1}}.
\end{equation}
For $a_{ij} = 0$ we have
\begin{equation}
    E_{i}E_{j} = E_{j}E_{i},
\end{equation}
\begin{equation}
    F_{i}F_{j} = F_{j}F_{i}.
\end{equation}
For $a_{ij} = -1$ we have
\begin{equation}
    E_{i}^{2}E_{j} - (q+q^{-1})E_{i}E_{j}E_{i} + E_{j}E_{i}^{2} = 0,
\end{equation}
\begin{equation}
    F_{i}^{2}F_{j} - (q+q^{-1})F_{i}F_{j}F_{i} + F_{j}F_{i}^{2} = 0,
\end{equation}
Coproduct is given by
\begin{equation}
\Delta E_{j} = E_{j}\otimes 1 + K_{j}^{-1}\otimes E_{j},
\end{equation}
\begin{equation}
\Delta F_{j} = 1\otimes F_{j} + F_{j}\otimes K_{j},
\end{equation}
\begin{equation}
\Delta K_{j} = K_{j}\otimes K_{j}.
\end{equation}

Non-compact quantum dilogarithm $G_{b}(z)$ is a special function introduced in \cite{F1} (see also \cite{F0}, \cite{FKV}, \cite{V}, \cite{Ka1}, \cite{KLSTS}, \cite{BT}). It is defined as follows
\begin{equation}
\log G_{b}(z) = \log\bar{\zeta}_{b} - \int\limits_{\mathbb{R}+\imath 0} \frac{dt}{t}\frac{e^{zt}}{(1-e^{bt})(1-e^{b^{-1}t})},
\end{equation}
where $Q = b+b^{-1}$ and $\zeta_{b} = e^{\frac{\pi\imath}{4} + \frac{\pi\imath(b^{2}+b^{-2})}{12}}$. Note, that $G_{b}(z)$ is closely related to the double sine function $S_{2}(z|\omega_{1},\omega_{2})$, see eq.(A.22) in \cite{KLSTS}.

Below we outline some properties of $G_{b}(z)$\\*
1. The function $G_{b}(z)$ has simple poles and zeros at the points
\begin{equation}
    z = -n_{1}b -n_{2}b^{-1},
\end{equation}
\begin{equation}
    z = Q +n_{1}b + n_{2}b^{-1},
\end{equation}
respectively, where $n_{1}$,$n_{2}$ are nonnegative integer numbers.\\*
2. $G_{b}(z)$ has the following asymptotic behavior:
\begin{equation}
 G_{b}(z) \sim
 \begin{cases} \bar{\zeta}_{b}, Im z \rightarrow +\infty ,\\ \zeta_{b} e^{\pi\imath z(z-Q)}, Im z \rightarrow -\infty . \end{cases}
\end{equation}
3. Functional equation:
\begin{equation}
G_{b}(z +b^{\pm 1}) = (1-e^{2\pi\imath b^{\pm 1}z})G_{b}(z).
\end{equation}
4. Reflection formula:
\begin{equation}
G_{b}(z)G_{b}(Q-z) = e^{\pi\imath z(z-Q)}.
\end{equation}\\*
5. 4-5 integral identity, \cite{V}:
\begin{equation}
\int d\tau
e^{-2\pi\gamma\tau}
\frac{G_{b}(\alpha+\imath\tau)G_{b}(\beta+\imath\tau)}{G_{b}(\alpha+\beta+\gamma+\imath\tau)G_{b}(Q+\imath\tau)} =
\frac{G_{b}(\alpha)G_{b}(\beta)G_{b}(\gamma)}{G_{b}(\alpha+\gamma)G_{b}(\beta+\gamma)}.
\end{equation}\\*
Define also the function $g_{b}(x)$ by
\begin{equation}
    g_{b}(x) = \frac{\bar{\zeta}_{b}}{G_{b}(\frac{Q}{2} +\frac{1}{2\pi\imath b}\log x)}.
\end{equation}
It has the following properties:\\*
1. $|g_{b}(x)| = 1$, if $x\in \mathbb{R}_{+}$. So if $A$ is a positive self-adjoint operator, then $g_{b}(A)$ is unitary. \\*
2. Fourier transform:
\begin{equation}
g_{b}(x) = \int d\tau x^{\imath b^{-1}\tau}e^{\pi Q\tau}G_{b}(-\imath\tau).
\end{equation}\\*
Let $U$, $V$ be positive self-adjoint operators satisfying the relation $UV = q^{2}VU$. Then the following non-commutative identities hold:\\*
3. Quantum exponential relation, \cite{F2}:
\begin{equation}
    g_{b}(U)g_{b}(V) = g_{b}(U+V).
\end{equation}
4. Quantum pentagon relation, \cite{Ka0}:
\begin{equation}
    g_{b}(V)g_{b}(U) = g_{b}(U)g_{b}(q^{-1}UV)g_{b}(V).
\end{equation}
5. Another useful relation, \cite{BT}:
\begin{equation}
    U + V = g_{b}(qU^{-1}V)Ug^{\ast}_{b}(qU^{-1}V),
\end{equation}
where the star means hermitian conjugation.

\section{Algebra of complex powers of generators of $U_{q}(\mathfrak{sl}(2))$}

Let $q = e^{\pi\imath b^{2}}$, $(b^{2}\in \mathbb{R}\setminus \mathbb{Q})$ and let $K = q^{H}$, $E$, $F$ be the generators of $U_{q}(\mathfrak{sl}(2))$ subjected to the relations
\begin{equation}
KE = q^{2}EK,
\end{equation}
\begin{equation}
KF = q^{-2}FK,
\end{equation}
\begin{equation}
EF-FE = \frac{K-K^{-1}}{q-q^{-1}}.
\end{equation}

Define the rescaled versions of generators $E$, $F$ by
\begin{equation}\label{rescaled E}
\mathcal{E} = -\imath(q-q^{-1})E,
\end{equation}
\begin{equation}
\mathcal{F} = -\imath(q-q^{-1})F.
\end{equation}

Let $\nu$ be a positive real number. There is a well-known representation of $U_{q}(\mathfrak{sl}(2))$ (see e.g.\cite{PT1}):
\begin{equation}
H = -2\imath b^{-1}u,
\end{equation}
\begin{equation}
K = q^{H} = e^{2\pi bu},
\end{equation}
\begin{equation}
\mathcal{E} =  q^{-\frac{1}{2}}e^{\pi b\nu +\pi bu}e^{-\imath b\partial_{u}} + q^{\frac{1}{2}}e^{-\pi b\nu -\pi bu}e^{-\imath b\partial_{u}},
\end{equation}
\begin{equation}
\mathcal{F} = q^{-\frac{1}{2}}e^{\pi b\nu-\pi bu}e^{\imath b\partial_{u}} +q^{\frac{1}{2}}e^{-\pi b\nu+\pi bu}e^{\imath b\partial_{u}}.
\end{equation}
This representation is a particular example of positive principal series representations of $U_{q}(\mathfrak{g})$ \cite{Ip2}.\\*
The following lemma was proven for a slightly different representation of $U_{q}(\mathfrak{sl}(2))$ in \cite{BT} and for the representation we use in this paper in \cite{Ip1}. It gives the expressions of the generators of $U_{q}(\mathfrak{sl}(2))$ in a form convenient for the definition of functions of them. It is based on the formula eq.(B.2) in \cite{BT}, stating that given positive self-adjoint operators $U$,$V$ satisfying $UV = q^{2}VU$, one can write
\begin{equation}
U+V = g_{b}(qU^{-1}V)U(g_{b}(qU^{-1}V))^{-1}.
\end{equation}
\begin{lem}
Let $\mathcal{E}$, $\mathcal{F}$ be the rescaled positive generators of $U_{q}(\mathfrak{sl}(2))$ defined above. They can be written in the following form:
\begin{equation}
\mathcal{E} = g_{b}(e^{-2\pi b\nu-2\pi bu})e^{\pi b\nu +\pi bu - \imath b\partial_{u}}g^{\ast}_{b}(e^{-2\pi b\nu-2\pi bu}),
\end{equation}
\begin{equation}
\mathcal{F} = g_{b}(e^{-2\pi b\nu+2\pi bu})e^{\pi b\nu-\pi bu+\imath b\partial_{u}}g^{\ast}_{b}(e^{-2\pi b\nu+2\pi bu}).
\end{equation}
\end{lem}
$\proof$
Let $U$, $V$ be positive self-adjoint operators such that
$$
UV = q^{2}VU.
$$
Then, \cite{BT}:
$$
U+V = g_{b}(qU^{-1}V)U(g_{b}(qU^{-1}V))^{-1}.
$$
For $\mathcal{E}$ we have $U = q^{-\frac{1}{2}}e^{\pi b\nu +\pi bu}e^{-\imath b\partial_{u}}$, $V = q^{\frac{1}{2}}e^{-\pi b\nu -\pi bu}e^{-\imath b\partial_{u}}$ and
$$
qU^{-1}V = q q^{\frac{1}{2}}e^{\imath b\partial_{u}}e^{-\pi b\nu-\pi bu}q^{\frac{1}{2}}e^{-\pi b\nu -\pi bu}e^{-\imath b\partial_{u}} = e^{-2\pi b\nu-2\pi bu},
$$
so
$$
\mathcal{E} = g_{b}(e^{-2\pi b\nu-2\pi bu})q^{-\frac{1}{2}}e^{\pi b\nu +\pi bu}e^{-\imath b\partial_{u}}(g_{b}(e^{-2\pi b\nu-2\pi bu}))^{-1}.
$$
For $\mathcal{F}$ we have $U = q^{-\frac{1}{2}}e^{\pi b\nu-\pi bu}e^{\imath b\partial_{u}}$, $V = q^{\frac{1}{2}}e^{-\pi b\nu+\pi bu}e^{\imath b\partial_{u}}$,
$$
qU^{-1}V = qq^{\frac{1}{2}}e^{-\imath b\partial_{u}}e^{-\pi b\nu+\pi bu}q^{\frac{1}{2}}e^{-\pi b\nu+\pi bu}e^{\imath b\partial_{u}} = e^{-2\pi b\nu +2\pi bu},
$$
$$
\mathcal{F} = g_{b}(e^{-2\pi b\nu+2\pi bu})q^{-\frac{1}{2}}e^{\pi b\nu-\pi bu}e^{\imath b\partial_{u}}(g_{b}(e^{-2\pi b\nu+2\pi bu}))^{-1}.
$$
$\Box$

Multiplication by $g_{b}(e^{-2\pi b\nu-2\pi bu})$ and $g_{b}(e^{-2\pi b\nu+2\pi bu})$ is unitary transformation, since $|g_{b}(x)| = 1$, for $x\in\mathbb{R}_{+}$. As a consequence, following \cite{BT}, eq.(3.15), eq.(3.21), one can define functions of generators $\mathcal{E}$ and $\mathcal{F}$ as follows:
\begin{de}
Let $\varphi(x)$ be a complex-valued function and let $K$, $\mathcal{E}$, $\mathcal{F}$ be $U_{q}(\mathfrak{sl}(2))$ generators in the positive principal series representation. The functions of these operators are defined as follows
\begin{equation}
    \varphi(K) = \varphi(e^{2\pi bu}),
\end{equation}
\begin{equation}\label{function of E}
\varphi(\mathcal{E}) = g_{b}(e^{-2\pi b\nu-2\pi bu})\varphi(e^{\pi b\nu +\pi bu - \imath b\partial_{u}})g^{\ast}_{b}(e^{-2\pi b\nu-2\pi bu}),
\end{equation}
\begin{equation}\label{function of F}
\varphi(\mathcal{F}) = g_{b}(e^{-2\pi b\nu+2\pi bu})\varphi(e^{\pi b\nu-\pi bu+\imath b\partial_{u}})g^{\ast}_{b}(e^{-2\pi b\nu+2\pi bu}).
\end{equation}
\end{de}
In particular, the powers of $\mathcal{E}$ and $\mathcal{F}$ are given by
\begin{equation}\label{Imaginary power E Uqsl(2)}
\mathcal{E}^{\imath s} = g_{b}(e^{-2\pi b\nu-2\pi bu})e^{\pi\imath bs\nu +\pi\imath bsu +  bs\partial_{u}}g^{\ast}_{b}(e^{-2\pi b\nu-2\pi bu}),
\end{equation}
\begin{equation}\label{Imaginary power F Uqsl(2)}
\mathcal{F}^{\imath t} = g_{b}(e^{-2\pi b\nu+2\pi bu})e^{\pi\imath bt\nu-\pi\imath btu- bt\partial_{u}}g^{\ast}_{b}(e^{-2\pi b\nu+2\pi bu}).
\end{equation}

The formulas for the powers in this particular representation were obtained in \cite{Ip1}.\\*
Define the arbitrary divided powers of $A$ by
\begin{equation}\label{complex devided power}
    A^{(\imath s)} = G_{b}(-\imath bs)A^{\imath s}.
\end{equation}

\begin{te}\cite{Su2}
The following generalized Kac's identity holds
\begin{equation}
\mathcal{E}^{(\imath s)}\mathcal{F}^{(\imath t)} = \int\limits_{\mathcal{C}} d\tau e^{\pi bQ\tau}\mathcal{F}^{(\imath t+\imath\tau)}K^{-\imath\tau}\frac{G_{b}(\imath b\tau)G_{b}(-bH+\imath b(s+t+\tau))}{G_{b}(-bH+\imath b(s+t+2\tau))}\mathcal{E}^{(\imath s+\imath\tau)},
\end{equation}
where the contour $\mathcal{C}$ goes slightly above the real axis but passes below the pole at $\tau = 0$.
\end{te}

\section{Algebra of complex powers of generators of $U_{q}(\mathfrak{sl}(3))$}

In this section we prove the Generalized Kac's identity in the case of positive principal series representations of $U_{q}(\mathfrak{sl}(3))$.

Let $q = e^{\pi\imath b^{2}}$, ($b^{2}\in \mathbb{R}\setminus \mathbb{Q}$). Let $a_{ij}$, ($i$,$j = 1$, $2$) be the Cartan matrix corresponding to $\mathfrak{sl}(3)$ Lie algebra, i.e. $a_{11} = a_{22} = 2$, $a_{12} = a_{21} = -1$. $U_{q}(\mathfrak{sl}(3))$ $(q = e^{\pi\imath b^{2}}$, $b^{2}\in \mathbb{R}\setminus \mathbb{Q})$ is defined by generators $E_{j}$, $F_{j}$, $K_{j} = q^{H_{j}}$, $1\le j \le 2$ and relations
\begin{equation}
    K_{i}K_{j} = K_{j}K_{i},
\end{equation}
\begin{equation}
    K_{i}E_{j} = q^{a_{ij}}E_{j}K_{i},
\end{equation}
\begin{equation}
    K_{i}F_{j} = q^{-a_{ij}}F_{j}K_{i},
\end{equation}
\begin{equation}
    E_{i}F_{j} - F_{j}E_{i} = \delta_{ij}\frac{K_{i} - K_{i}^{-1}}{q-q^{-1}}.
\end{equation}
For $i\ne j$ we have
\begin{equation}
    E_{i}^{2}E_{j} - (q+q^{-1})E_{i}E_{j}E_{i} + E_{j}E_{i}^{2} = 0,
\end{equation}
\begin{equation}
    F_{i}^{2}F_{j} - (q+q^{-1})F_{i}F_{j}F_{i} + F_{j}F_{i}^{2} = 0.
\end{equation}

The general construction of the positive principal series representations of $U_{q}(\mathfrak{g})$ in the simply-laced case using  Lusztig's data was given in \cite{Ip2}. Let $w_{0}$ be the longest element of the Weyl group. There are different realizations of positive principal series representations corresponding to each reduced expressions of $w_{0}$. In the case of $U_{q}(\mathfrak{sl}(3))$ there are two options $w_{0} = s_{1}s_{2}s_{1}$ and $w_{0} = s_{2}s_{1}s_{2}$. In the following we give the explicit formulas for both of this cases.

Let $\mathcal{E}_{j} = -\imath(q-q^{-1})E_{j}$, $\mathcal{F}_{j} = -\imath(q-q^{-1})F_{j}$, $j = 1$, $2$ be the rescaled versions of $U_{q}(\mathfrak{sl}(3))$ generators.

\begin{prop} \cite{Ip2}.
Let $K_{j}$, $\mathcal{E}_{j}$, $\mathcal{F}_{j}$ be the rescaled generators of $U_{q}(\mathfrak{sl}(3))$. Let $w_{0} = s_{1}s_{2}s_{1}$ be reduced expression of the longest Weyl element. Let $\nu_{1}$,$\nu_{2}$ be positive real numbers. The positive principal series representation of $U_{q}(\mathfrak{sl}(3))$ corresponding to these data is given by:
\begin{equation}
K_{1} = e^{-2\pi b\nu_{1}+2\pi bu-\pi bv+2\pi bw},
\end{equation}
\begin{equation}
K_{2} = e^{-2\pi b\nu_{2}-\pi bu+2\pi bv-\pi bw},
\end{equation}
\begin{equation}
\mathcal{E}_{1} = e^{\pi bw-\imath b\partial_{w}} + e^{-\pi bw-\imath b\partial_{w}},
\end{equation}
\begin{equation}
\mathcal{E}_{2} = e^{\pi bv-\pi bw-\imath b\partial_{v}}+ e^{\pi bu -\imath b\partial_{u} -\imath b\partial_{v}+\imath b\partial_{w}} +
e^{-\pi bu-\imath b\partial_{u}-\imath b\partial_{v}+\imath b\partial_{w}} +
e^{-\pi bv+\pi bw-\imath b\partial_{v}},
\end{equation}
\begin{equation}
\mathcal{F}_{1} =
e^{2\pi b\nu_{1}-2\pi bu+\pi bv-\pi bw+\imath b\partial_{w}} + e^{2\pi b\nu_{1}-\pi bu+\imath b\partial_{u}} +
e^{-2\pi b\nu_{1}+\pi bu+\imath b\partial_{u}} + e^{-2\pi b\nu_{1}+2\pi bu-\pi bv+\pi bw+\imath b\partial_{w}},
\end{equation}
\begin{equation}
\mathcal{F}_{2} = e^{2\pi b\nu_{2}+\pi bu-\pi bv +\imath b\partial_{v}} + e^{-2\pi b\nu_{2}-\pi bu +\pi bv+\imath b\partial_{v}}.
\end{equation}
\end{prop}

Similar to $U_{q}(\mathfrak{sl}(2))$ case using eq.(B.2), \cite{BT} we represent the generators in a form convenient for the definition of functions of them
\begin{lem}
Let $\mathcal{E}_{i}$, $\mathcal{F}_{i}$, $(i = 1,2)$ be the  generators of $U_{q}(\mathfrak{sl}(3))$ in the positive principal series representation corresponding to the reduced expression $w_{0} = s_{1}s_{2}s_{1}$. They can be represented in the following form
\begin{equation}
\mathcal{E}_{1} = g_{b}(e^{-2\pi bw})e^{\pi bw-\imath b\partial_{w}}g_{b}^{\ast}(e^{-2\pi bw}),
\end{equation}
\begin{multline}
\mathcal{E}_{2} =
g_{b}(e^{\pi bu-\pi bv +\pi bw-\imath b\partial_{u}+\imath b\partial_{w}})g_{b}(e^{-\pi bu -\pi bv+\pi bw-\imath b\partial_{u}+\imath b\partial_{w}})g_{b}(e^{-2\pi bv+2\pi bw})\times \\
e^{\pi bv-\pi bw-\imath b\partial_{v}}\times      \\
g^{\ast}_{b}(e^{-2\pi bv+2\pi bw})g^{\ast}_{b}(e^{-\pi bu -\pi bv+\pi bw-\imath b\partial_{u}+\imath b\partial_{w}})
g^{\ast}_{b}(e^{\pi bu-\pi bv +\pi bw-\imath b\partial_{u}+\imath b\partial_{w}}),
\end{multline}
\begin{multline}
\mathcal{F}_{1} = g_{b}(e^{\pi bu-\pi bv+\pi bw+\imath b\partial_{u}-\imath b\partial_{w}})
g_{b}(e^{-4\pi b\nu_{1}+3\pi bu-\pi bv+\pi bw+\imath b\partial_{u}-\imath b\partial_{w}})
g_{b}(e^{-4\pi b\nu_{1}+4\pi bu-2\pi bv+2\pi bw})\times \\
e^{2\pi b\nu_{1}-2\pi bu+\pi bv-\pi bw+\imath b\partial_{w}}\times \\
g_{b}^{\ast}(e^{-4\pi b\nu_{1}+4\pi bu-2\pi bv+2\pi bw})
g_{b}^{\ast}(e^{-4\pi b\nu_{1}+3\pi bu-\pi bv+\pi bw+\imath b\partial_{u}-\imath b\partial_{w}})
g_{b}^{\ast}(e^{\pi bu-\pi bv+\pi bw+\imath b\partial_{u}-\imath b\partial_{w}}),
\end{multline}
\begin{equation}
\mathcal{F}_{2} = g_{b}(e^{-4\pi b\nu_{2}-2\pi bu+2\pi bv})e^{2\pi b\nu_{2}+\pi bu-\pi bv+\imath b\partial_{v}}
g^{\ast}_{b}(e^{-4\pi b\nu_{2}-2\pi bu+2\pi bv}).
\end{equation}
\end{lem}
$\proof$

Let $q=e^{\pi\imath b^{2}}$ and let $U$, $V$ be positive essentially self-adjoint operators subjected to the relation $UV=q^{2}VU$. We will need the following identity,\cite{BT}:
$$
U+V = g_{b}(qU^{-1}V)Ug^{\ast}_{b}(qU^{-1}V),
$$
and the quantum exponential relation, \cite{F2}:
$$
g_{b}(U+V) = g_{b}(U)g_{b}(V).
$$
Let us start with $\mathcal{E}_{1}$. It has the following form
$$
\mathcal{E}_{1} = U+V,
$$
where $U = e^{\pi bw-\imath b\partial_{w}}$, $V = e^{-\pi bw-\imath b\partial_{w}}$. Using the identities $e^{A}e^{B} = e^{\frac{[A,B]}{2}}e^{A+B}$ and $e^{A}e^{B} = e^{[A,B]}e^{B}e^{A}$ in the case when the commutator $[A,B]$ commutes with both $A$ and $B$, and also the identity $[x,\partial_{x}] = -1$, one checks that
$$
UV = e^{\pi bw-\imath b\partial_{w}}e^{-\pi bw-\imath b\partial_{w}} =
e^{[\pi bw-\imath b\partial_{w},-\pi bw-\imath b\partial_{w}]}e^{-\pi bw-\imath b\partial_{w}}e^{\pi bw-\imath b\partial_{w}} =
$$
$$
e^{2\pi\imath b^{2}}e^{-\pi bw-\imath b\partial_{w}}e^{\pi bw-\imath b\partial_{w}} = q^{2}VU,
$$
and
$$
qU^{-1}V = e^{\pi\imath b^{2}}e^{-\pi bw+\imath b\partial_{w}}e^{-\pi bw-\imath b\partial_{w}}=
e^{\pi\imath b^{2}}e^{\frac{1}{2}[-\pi bw+\imath b\partial_{w},-\pi bw-\imath b\partial_{w}]} = e^{-2\pi bw},
$$
so we have
$$
\mathcal{E}_{1} = g_{b}(qU^{-1}V)Ug^{\ast}_{b}(qU^{-1}V) = g_{b}(e^{-2\pi bw})e^{\pi bw-\imath b\partial_{w}}g_{b}^{\ast}(e^{-2\pi bw}).
$$
For $\mathcal{E}_{2}$ we have
$$
\mathcal{E}_{2} = U_{1} + U_{2} + U_{3} + U_{4},
$$
where $U_{1} = e^{\pi bv-\pi bw-\imath b\partial_{v}}$ $U_{2} = e^{\pi bu -\imath b\partial_{u} -\imath b\partial_{v}+\imath b\partial_{w}}$, $U_{3} = e^{-\pi bu-\imath b\partial_{u}-\imath b\partial_{v}+\imath b\partial_{w}}$,  $U_{4} =  e^{-\pi bv+\pi bw-\imath b\partial_{v}}$.
These operators satisfy the relations
$$
U_{i}U_{j} = q^{2}U_{j}U_{i},
$$
if $i< j$. We obtain
$$
\mathcal{E}_{2} = g_{b}(qU_{1}^{-1}(U_{2}+U_{3}+U_{4}))U_{1}g^{\ast}_{b}(qU_{1}^{-1}(U_{2}+U_{3}+U_{4})) =
$$
$$
g_{b}(qU_{1}^{-1}U_{2})g_{b}(qU_{1}^{-1}U_{3})g_{b}(qU_{1}^{-1}U_{4})U_{1}g^{\ast}_{b}(qU_{1}^{-1}U_{4})g^{\ast}_{b}(qU_{1}^{-1}U_{3})g^{\ast}_{b}(qU_{1}^{-1}U_{2}),
$$
where in the second equality we have used the quantum exponential relation, provided that operators $(qU_{1}^{-1}U_{i})$ are positive and satisfy $(qU_{1}^{-1}U_{i})(qU_{1}^{-1}U_{j}) = q^{2}(qU_{1}^{-1}U_{j})(qU_{1}^{-1}U_{i})$ for $1<i<j$ .
$$
qU_{1}^{-1}U_{2} = e^{\pi\imath b^{2}}e^{-\pi bv+\pi bw+\imath b\partial_{v}}e^{\pi bu -\imath b\partial_{u} -\imath b\partial_{v}+\imath b\partial_{w}} =
$$
$$
e^{\pi\imath b^{2}}e^{\frac{1}{2}[-\pi bv+\pi bw+\imath b\partial_{v},\pi bu -\imath b\partial_{u} -\imath b\partial_{v}+\imath b\partial_{w}]}e^{\pi bu-\pi bv+\pi bw-\imath b\partial_{u}+\imath b\partial_{w}} =
$$
$$
e^{\pi\imath b^{2}}e^{-\pi\imath b^{2}}e^{\pi bu-\pi bv+\pi bw-\imath b\partial_{u}+\imath b\partial_{w}} =
e^{\pi bu-\pi bv+\pi bw-\imath b\partial_{u}+\imath b\partial_{w}}.
$$
Analogously we obtain
$$
qU_{1}^{-1}U_{3} = e^{-\pi bu-\pi bv+\pi bw-\imath b\partial_{u}+\imath b\partial_{w}},
$$
$$
qU_{1}^{-1}U_{4}= e^{-2\pi bv+2\pi bw}.
$$
Substituting these expressions into the formula for $\mathcal{E}_{2}$ we obtain
$$
\mathcal{E}_{2} = g_{b}(e^{\pi bu-\pi bv +\pi bw-\imath b\partial_{u}+\imath b\partial_{w}})
g_{b}(e^{-\pi bu -\pi bv+\pi bw-\imath b\partial_{u}+\imath b\partial_{w}})
g_{b}(e^{-2\pi bv+2\pi bw})e^{\pi bv-\pi bw-\imath b\partial_{v}}\times
$$
$$
g^{\ast}_{b}(e^{-2\pi bv+2\pi bw})g^{\ast}_{b}(e^{-\pi bu -\pi bv+\pi bw-\imath b\partial_{u}+\imath b\partial_{w}})
g^{\ast}_{b}(e^{\pi bu-\pi bv +\pi bw-\imath b\partial_{u}+\imath b\partial_{w}}).
$$

For the generator $\mathcal{F}_{1}$ we have
$$
\mathcal{F}_{1} = U_{1} + U_{2} + U_{3} + U_{4},
$$
where we used the notations
$$
U_{1} = e^{2\pi b\nu_{1}-2\pi bu+\pi bv-\pi bw+\imath b\partial_{w}},$$
$$
U_{2} = e^{2\pi b\nu_{1}-\pi bu+\imath b\partial_{u}},
$$
$$
U_{3} = e^{-2\pi b\nu_{1}+\pi bu+\imath b\partial_{u}},
$$
$$
U_{4} = e^{-2\pi b\nu_{1}+2\pi bu-\pi bv+\pi bw+\imath b\partial_{w}}.
$$
Again, for $i<j$ we have the relations
$$
U_{i}U_{j} = q^{2}U_{j}U_{i},
$$
and for $1<i<j$
$$
(qU^{-1}_{1}U_{i})(qU^{-1}_{1}U_{j}) = q^{2}(qU^{-1}_{1}U_{j})(qU^{-1}_{1}U_{i}).
$$
Explicit expressions for the operators $qU^{-1}_{1}U_{i}$ are given by
$$
qU_{1}^{-1}U_{2} = e^{\pi bu-\pi bv+\pi bw+\imath b\partial_{u}-\imath b\partial_{w}},
$$
$$
qU_{1}^{-1}U_{3} = e^{-4\pi b\nu_{1}+3\pi bu-\pi bv+\pi bw+\imath b\partial_{u}-\imath b\partial_{w}},
$$
$$
qU_{1}^{-1}U_{4} = e^{-4\pi b\nu_{1}+4\pi bu-2\pi bv+2\pi bw}.
$$
Using the formulas $U+V = g_{b}(qU^{-1}V)Ug^{\ast}_{b}(qU^{-1}V)$ and $g_{b}(U+V) = g_{b}(U)g_{b}(V)$ for positive operators satisfying the relation $UV = q^{2}VU$ we obtain
$$
\mathcal{F}_{1} = g_{b}(qU_{1}^{-1}(U_{2}+U_{3}+U_{4}))U_{1}g^{\ast}_{b}(qU_{1}^{-1}(U_{2}+U_{3}+U_{4})) =
$$
$$
g_{b}(qU_{1}^{-1}U_{2})g_{b}(qU_{1}^{-1}U_{3})g_{b}(qU_{1}^{-1}U_{4})U_{1}g^{\ast}_{b}(qU_{1}^{-1}U_{4})g^{\ast}_{b}(qU_{1}^{-1}U_{3})g^{\ast}_{b}(qU_{1}^{-1}U_{2}) =
$$
$$
g_{b}(e^{\pi bu-\pi bv+\pi bw+\imath b\partial_{u}-\imath b\partial_{w}})
g_{b}(e^{-4\pi b\nu_{1}+3\pi bu-\pi bv+\pi bw+\imath b\partial_{u}-\imath b\partial_{w}})
g_{b}(e^{-4\pi b\nu_{1}+4\pi bu-2\pi bv+2\pi bw})\times
$$
$$
e^{2\pi b\nu_{1}-2\pi bu+\pi bv-\pi bw+\imath b\partial_{w}}\times
$$
$$
g_{b}^{\ast}(e^{-4\pi b\nu_{1}+4\pi bu-2\pi bv+2\pi bw})
g_{b}^{\ast}(e^{-4\pi b\nu_{1}+3\pi bu-\pi bv+\pi bw+\imath b\partial_{u}-\imath b\partial_{w}})
g_{b}^{\ast}(e^{\pi bu-\pi bv+\pi bw+\imath b\partial_{u}-\imath b\partial_{w}}).
$$

Generator $\mathcal{F}_{2}$.
$$
\mathcal{F}_{2} = A_{1} + A_{2},
$$
where
$$
A_{1} = e^{2\pi b\nu_{2}+\pi bu-\pi bv +\imath b\partial_{v}},
$$
$$
A_{2} = e^{-2\pi b\nu_{2}-\pi bu +\pi bv+\imath b\partial_{v}}.
$$
Then
$$
qA_{1}^{-1}A_{2} = e^{-4\pi b\nu_{2}-2\pi bu+2\pi bv},
$$
$$
A_{1}A_{2} = q^{2}A_{2}A_{1},
$$
and using the identity
$$
A_{1} + A_{2} = g_{b}(qA_{1}^{-1}A_{2})A_{1}g_{b}^{\ast}(qA_{1}^{-1}A_{2}),
$$
we obtain
$$
\mathcal{F}_{2} = g_{b}(e^{-4\pi b\nu_{2}-2\pi bu+2\pi bv})e^{2\pi b\nu_{2}+\pi bu-\pi bv+\imath b\partial_{v}}
g^{\ast}_{b}(e^{-4\pi b\nu_{2}-2\pi bu+2\pi bv}).
$$

$\Box$

Similar to eq.(3.15), eq.(3.21) \cite{BT}, we define the functions of the operators in the following way:
\begin{de}
Let $\varphi(x)$ be a complex-valued function of one variable and let $K_{i}$, $\mathcal{E}_{i}$, $\mathcal{F}_{i}$, $(i = 1,2)$ be the  generators of $U_{q}(\mathfrak{sl}(3))$ in the positive principal series representation corresponding to the reduced expression $w_{0} = s_{1}s_{2}s_{1}$.
\begin{equation}
\varphi(K_{1}) = \varphi(e^{-2\pi b\nu_{1}+2\pi bu-\pi bv+2\pi bw}),
\end{equation}

\begin{equation}
\varphi(K_{2}) = \varphi(e^{-2\pi b\nu_{2}-\pi bu+2\pi bv-\pi bw}),
\end{equation}

\begin{equation}
\varphi(\mathcal{E}_{1}) = g_{b}(e^{-2\pi bw})\varphi(e^{\pi bw-\imath b\partial_{w}})g_{b}^{\ast}(e^{-2\pi bw}),
\end{equation}

\begin{multline}
\varphi(\mathcal{E}_{2}) = g_{b}(e^{\pi bu-\pi bv +\pi bw-\imath b\partial_{u}+\imath b\partial_{w}})
g_{b}(e^{-\pi bu -\pi bv+\pi bw-\imath b\partial_{u}+\imath b\partial_{w}})
g_{b}(e^{-2\pi bv+2\pi bw})\times  \\
\varphi(e^{\pi bv-\pi bw-\imath b\partial_{v}})\times \\
g^{\ast}_{b}(e^{-2\pi bv+2\pi bw})g^{\ast}_{b}(e^{-\pi bu -\pi bv+\pi bw-\imath b\partial_{u}+\imath b\partial_{w}})
g^{\ast}_{b}(e^{\pi bu-\pi bv +\pi bw-\imath b\partial_{u}+\imath b\partial_{w}}),
\end{multline}

\begin{multline}
\varphi(\mathcal{F}_{1}) = g_{b}(e^{\pi bu-\pi bv+\pi bw+\imath b\partial_{u}-\imath b\partial_{w}})
g_{b}(e^{-4\pi b\nu_{1}+3\pi bu-\pi bv+\pi bw+\imath b\partial_{u}-\imath b\partial_{w}})
g_{b}(e^{-4\pi b\nu_{1}+4\pi bu-2\pi bv+2\pi bw})\times  \\
\varphi(e^{2\pi b\nu_{1}-2\pi bu+\pi bv-\pi bw+\imath b\partial_{w}})\times \\
g_{b}^{\ast}(e^{-4\pi b\nu_{1}+4\pi bu-2\pi bv+2\pi bw})
g_{b}^{\ast}(e^{-4\pi b\nu_{1}+3\pi bu-\pi bv+\pi bw+\imath b\partial_{u}-\imath b\partial_{w}})
g_{b}^{\ast}(e^{\pi bu-\pi bv+\pi bw+\imath b\partial_{u}-\imath b\partial_{w}}),
\end{multline}

\begin{equation}
\varphi(\mathcal{F}_{2}) = g_{b}(e^{-4\pi b\nu_{2}-2\pi bu+2\pi bv})\varphi(e^{2\pi b\nu_{2}+\pi bu-\pi bv+\imath b\partial_{v}})
g^{\ast}_{b}(e^{-4\pi b\nu_{2}-2\pi bu+2\pi bv}).
\end{equation}
\end{de}
Choosing in the definition the function to be $\varphi(x) = x^{\imath s}$ we obtain the expressions for arbitrary powers of generators.

In the following we repeat the same steps for representations corresponding to another choice of reduced expression $w_{0} = s_{2}s_{1}s_{2}$.
\begin{prop}\cite{Ip2}.
The positive principal series representation of $U_{q}(\mathfrak{sl}(3))$ corresponding to the reduced expression of the longest Weyl element $w_{0} = s_{2}s_{1}s_{2}$ and positive real parameters $\nu_{1}$,$\nu_{2}$ is given by
\begin{equation}
K_{1} = e^{-2\pi b\nu_{1}-\pi bu+2\pi bv-\pi bw},
\end{equation}

\begin{equation}
K_{2} = e^{-2\pi b\nu_{2}+2\pi bu-\pi bv+2\pi bw},
\end{equation}

\begin{equation}
\mathcal{E}_{1} = e^{\pi bv-\pi bw-\imath b\partial_{v}}+ e^{\pi bu -\imath b\partial_{u} -\imath b\partial_{v}+\imath b\partial_{w}} +
e^{-\pi bu-\imath b\partial_{u}-\imath b\partial_{v}+\imath b\partial_{w}} +
e^{-\pi bv+\pi bw-\imath b\partial_{v}},
\end{equation}

\begin{equation}
\mathcal{E}_{2} = e^{\pi bw-\imath b\partial_{w}} + e^{-\pi bw-\imath b\partial_{w}},
\end{equation}

\begin{equation}
\mathcal{F}_{1} = e^{2\pi b\nu_{1}+\pi bu-\pi bv +\imath b\partial_{v}} + e^{-2\pi b\nu_{1}-\pi bu +\pi bv+\imath b\partial_{v}},
\end{equation}

\begin{equation}
\mathcal{F}_{2} =
e^{2\pi b\nu_{2}-2\pi bu+\pi bv-\pi bw+\imath b\partial_{w}} + e^{2\pi b\nu_{2}-\pi bu+\imath b\partial_{u}} +
e^{-2\pi b\nu_{2}+\pi bu+\imath b\partial_{u}} + e^{-2\pi b\nu_{2}+2\pi bu-\pi bv+\pi bw+\imath b\partial_{w}}.
\end{equation}
\end{prop}

\begin{lem}
Let $\mathcal{E}_{i}$, $\mathcal{F}_{i}$, $(i = 1,2)$ be the  generators of $U_{q}(\mathfrak{sl}(3))$ in the positive principal series representation corresponding to the reduced expression $w_{0} = s_{2}s_{1}s_{2}$. They can be represented in the following form

\begin{multline}
\mathcal{E}_{1} = g_{b}(e^{\pi bu-\pi bv +\pi bw-\imath b\partial_{u}+\imath b\partial_{w}})
g_{b}(e^{-\pi bu -\pi bv+\pi bw-\imath b\partial_{u}+\imath b\partial_{w}})
g_{b}(e^{-2\pi bv+2\pi bw})\times  \\
e^{\pi bv-\pi bw-\imath b\partial_{v}}\times  \\
g^{\ast}_{b}(e^{-2\pi bv+2\pi bw})g^{\ast}_{b}(e^{-\pi bu -\pi bv+\pi bw-\imath b\partial_{u}+\imath b\partial_{w}})
g^{\ast}_{b}(e^{\pi bu-\pi bv +\pi bw-\imath b\partial_{u}+\imath b\partial_{w}}),
\end{multline}

\begin{equation}
\mathcal{E}_{2} = g_{b}(e^{-2\pi bw})e^{\pi bw-\imath b\partial_{w}}g_{b}^{\ast}(e^{-2\pi bw}),
\end{equation}

\begin{equation}
\mathcal{F}_{1} = g_{b}(e^{-4\pi b\nu_{1}-2\pi bu+2\pi bv})e^{2\pi b\nu_{1} +\pi bu-\pi bv+\imath b\partial_{v}}
g^{\ast}_{b}(e^{-4\pi b\nu_{1}-2\pi bu+2\pi bv}),
\end{equation}

\begin{multline}
\mathcal{F}_{2} = g_{b}(e^{\pi bu-\pi bv+\pi bw+\imath b\partial_{u}-\imath b\partial_{w}})
g_{b}(e^{-4\pi b\nu_{2}+3\pi bu-\pi bv+\pi bw+\imath b\partial_{u}-\imath b\partial_{w}})
g_{b}(e^{-4\pi b\nu_{2}+4\pi bu-2\pi bv+2\pi bw})\times  \\
e^{2\pi b\nu_{2}-2\pi bu+\pi bv-\pi bw+\imath b\partial_{w}}\times \\
g_{b}^{\ast}(e^{-4\pi b\nu_{2}+4\pi bu-2\pi bv+2\pi bw})
g_{b}^{\ast}(e^{-4\pi b\nu_{2}+3\pi bu-\pi bv+\pi bw+\imath b\partial_{u}-\imath b\partial_{w}})
g_{b}^{\ast}(e^{\pi bu-\pi bv+\pi bw+\imath b\partial_{u}-\imath b\partial_{w}}).
\end{multline}
\end{lem}

\begin{de}
Let $\varphi(x)$ be a complex-valued function of one variable and let $K_{i}$, $\mathcal{E}_{i}$, $\mathcal{F}_{i}$, $(i = 1,2)$ be the  generators of $U_{q}(\mathfrak{sl}(3))$ in the positive principal series representation corresponding to the reduced expression $w_{0} = s_{2}s_{1}s_{2}$.The functions of the generators are defined as follows:

\begin{equation}
\varphi(K_{1}) = \varphi(e^{-2\pi b\nu_{1}-\pi bu+2\pi bv-\pi bw}),
\end{equation}

\begin{equation}
\varphi(K_{2}) = \varphi(e^{-2\pi b\nu_{2}+2\pi bu-\pi bv+2\pi bw}),
\end{equation}

\begin{multline}
\varphi(\mathcal{E}_{1}) = g_{b}(e^{\pi bu-\pi bv +\pi bw-\imath b\partial_{u}+\imath b\partial_{w}})
g_{b}(e^{-\pi bu -\pi bv+\pi bw-\imath b\partial_{u}+\imath b\partial_{w}})
g_{b}(e^{-2\pi bv+2\pi bw})\times \\
\varphi(e^{\pi bv-\pi bw-\imath b\partial_{v}})\times \\
g^{\ast}_{b}(e^{-2\pi bv+2\pi bw})g^{\ast}_{b}(e^{-\pi bu -\pi bv+\pi bw-\imath b\partial_{u}+\imath b\partial_{w}})
g^{\ast}_{b}(e^{\pi bu-\pi bv +\pi bw-\imath b\partial_{u}+\imath b\partial_{w}}),
\end{multline}

\begin{equation}
\varphi(\mathcal{E}_{2}) = g_{b}(e^{-2\pi bw})\varphi(e^{\pi bw-\imath b\partial_{w}})g_{b}^{\ast}(e^{-2\pi bw}),
\end{equation}

\begin{equation}
\varphi(\mathcal{F}_{1}) = g_{b}(e^{-4\pi b\nu_{1}-2\pi bu+2\pi bv})\varphi(e^{2\pi b\nu_{1} +\pi bu-\pi bv+\imath b\partial_{v}})
g^{\ast}_{b}(e^{-4\pi b\nu_{1}-2\pi bu+2\pi bv}),
\end{equation}

\begin{multline}
\varphi(\mathcal{F}_{2}) = g_{b}(e^{\pi bu-\pi bv+\pi bw+\imath b\partial_{u}-\imath b\partial_{w}})
g_{b}(e^{-4\pi b\nu_{2}+3\pi bu-\pi bv+\pi bw+\imath b\partial_{u}-\imath b\partial_{w}})
g_{b}(e^{-4\pi b\nu_{2}+4\pi bu-2\pi bv+2\pi bw})\times  \\
\varphi(e^{2\pi b\nu_{2}-2\pi bu+\pi bv-\pi bw+\imath b\partial_{w}})\times \\
g_{b}^{\ast}(e^{-4\pi b\nu_{2}+4\pi bu-2\pi bv+2\pi bw})
g_{b}^{\ast}(e^{-4\pi b\nu_{2}+3\pi bu-\pi bv+\pi bw+\imath b\partial_{u}-\imath b\partial_{w}})
g_{b}^{\ast}(e^{\pi bu-\pi bv+\pi bw+\imath b\partial_{u}-\imath b\partial_{w}}).
\end{multline}

\end{de}
Assuming $\varphi(x) = x^{\imath s}$, we obtain the expressions for arbitrary powers of generators.\\*
Recall the definition of the divided powers of $A$
\begin{equation}
    A^{(\imath s)} = G_{b}(-\imath bs)A^{\imath s}.
\end{equation}

Now, after we have defined arbitrary devided powers of the generators of $U_{q}(\mathfrak{sl}(3))$ in the positive principal series representations corresponding to both reduced expressions of the Weyl element, we can state the main theorem
\begin{te}
Let $q = e^{\pi\imath b^{2}}$, $(b^{2}\in \mathbb{R}\setminus \mathbb{Q})$ and  let $K_{j} = q^{H_{j}}$, $\mathcal{E}_{j} = -\imath (q-q^{-1})E_{j}$, $\mathcal{F}_{j} = -\imath (q-q^{-1})F_{j}$, $1\le j\le 2$ be $U_{q}(\mathfrak{sl}(3))$ generators in the positive principal series representation corresponding to any reduced expression of the Weyl element. Then the following generalized Kac's identity holds:
\bigskip
\begin{equation}\begin{split}
\mathcal{E}_{j}^{(\imath s)}\mathcal{F}_{j}^{(\imath t)} = \int\limits_{\mathcal{C}} d\tau e^{\pi bQ\tau}\mathcal{F}_{j}^{(\imath t+\imath\tau)}K_{j}^{-\imath \tau}
\frac{G_{b}(\imath b\tau)G_{b}(-bH_{j} + \imath b(s+t+\tau))}{G_{b}(-bH_{j}+\imath b(s+t+2\tau))}\mathcal{E}_{j}^{(\imath s + \imath \tau)},
\end{split}\end{equation}
\bigskip
where the contour $\mathcal{C}$ goes slightly above the real axis but passes below the pole at $\tau = 0$.
\end{te}
$\proof$
The proof follows from the results stated in Proposition 4.3, Lemma 4.3, Proposition 4.4, Corollary 4.1, Corollary 4.2, Corollary 4.3.

 The next statement (Theorem 5.7 in \cite{Ip2}) establishes the unitary equivalence of positive principal series representations corresponding to different expressions of the longest Weyl element. We give here another proof for the case of $U_{q}(\mathfrak{sl}(3))$ of this result which allows explicitly illustrate its validness for arbitrary functions of generators. In this proof the pentagon identity \cite{Ka0} is extensively used which states that for positive self-adjoint operators $U$,$V$ satisfying the relation $UV = q^{2}VU$ we have
\begin{equation}
    g_{b}(V)g_{b}(U) = g_{b}(U)g_{b}(q^{-1}UV)g_{b}(V).
\end{equation}

\begin{prop}
Let $X_{s_{1}s_{2}s_{1}}$ be any generator of $U_{q}(\mathfrak{sl}(3))$ in the positive principal series representation corresponding to the reduced expression $w_{0} = s_{1}s_{2}s_{1}$. Let $X_{s_{2}s_{1}s_{2}}$ be the same generator in the positive principal series representation corresponding to the reduced expression $w_{0} = s_{2}s_{1}s_{2}$ and let $\varphi(x)$ be a complex-valued function.
The unitary transformation defined by
\begin{equation}
U = (uv)(vw)e^{-w\partial_{u}}e^{w\partial_{v}}g_{b}(e^{-\pi bu+\pi bv-\pi bw-\imath b\partial_{u}+\imath b\partial_{w}})
g^{\ast}_{b}(e^{\pi bu-\pi bv+\pi bw-\imath b\partial_{u}+\imath b\partial_{w}}).
\end{equation}
relates the functions of generators in these two representations by
$$
\varphi(X_{s_{2}s_{1}s_{2}}) = U\varphi(X_{s_{1}s_{2}s_{1}})U^{\ast}.
$$
\end{prop}
$\proof$
Let $\mathcal{E}_{1}$ be the generator in $s_{1}s_{2}s_{1}$ representation. Its unitary transform is given by
$$
U\varphi(\mathcal{E}_{1})U^{\ast} =
(uv)(vw)e^{-w\partial_{u}}e^{w\partial_{v}}g_{b}(e^{-\pi bu+\pi bv-\pi bw-\imath b\partial_{u}+\imath b\partial_{w}})
g^{\ast}_{b}(e^{\pi bu-\pi bv+\pi bw-\imath b\partial_{u}+\imath b\partial_{w}})\times
$$
$$
g_{b}(e^{-2\pi bw})\varphi(e^{\pi bw-\imath b\partial_{w}})(h.c.) =
$$
$$
(uv)(vw)e^{-w\partial_{u}}e^{w\partial_{v}}g^{\ast}_{b}(e^{\pi bu-\pi bv+\pi bw-\imath b\partial_{u}+\imath b\partial_{w}})g_{b}(e^{-\pi bu+\pi bv-\pi bw-\imath b\partial_{u}+\imath b\partial_{w}})
\times
$$
$$
g_{b}(e^{-2\pi bw})\varphi(e^{\pi bw-\imath b\partial_{w}})(h.c.).
$$
We have used the commutation of the factors $g^{\ast}_{b}(e^{\pi bu-\pi bv+\pi bw-\imath b\partial_{u}+\imath b\partial_{w}})$ and $g_{b}(e^{-\pi bu+\pi bv-\pi bw-\imath b\partial_{u}+\imath b\partial_{w}})$. Now, according to the identity $g_{b}(x)g_{b}(\frac{1}{x}) = e^{\frac{\pi\imath}{4\pi^{2}b^{2}}\log^{2}x}$ we have
$$
g^{\ast}_{b}(e^{\pi bu-\pi bv+\pi bw-\imath b\partial_{u}+\imath b\partial_{w}}) = e^{-\frac{\pi\imath}{4\pi^2 b^{2}}(\pi bu-\pi bv+\pi bw-\imath b\partial_{u}+\imath b\partial_{w})^{2}}
g_{b}(e^{-\pi bu+\pi bv-\pi bw+\imath b\partial_{u}-\imath b\partial_{w}}).
$$
Let
$$
A_{1} = e^{-2\pi bw},
$$
$$
A_{2} = e^{-\pi bu+\pi bv-\pi bw-\imath b\partial_{u}+\imath b\partial_{w}},
$$
Then
$$
q^{-1}A_{1}A_{2} = e^{-\pi bu+\pi bv-3\pi bw-\imath b\partial_{u}+\imath b\partial_{w}}
$$
and $A_{1}A_{2} = q^{2}A_{2}A_{1}$. Using the pentagon identity, \cite{Ka0}:
$$
g_{b}(A_{2})g_{b}(A_{1}) = g_{b}(A_{1})g_{b}(q^{-1}A_{1}A_{2})g_{b}(A_{2}),
$$
we have
$$
g_{b}(e^{-\pi bu+\pi bv-\pi bw-\imath b\partial_{u}+\imath b\partial_{w}})g_{b}(e^{-2\pi bw}) = g_{b}(e^{-2\pi bw})g_{b}(e^{-\pi bu+\pi bv-3\pi bw-\imath b\partial_{u}+\imath b\partial_{w}})g_{b}(e^{-\pi bu+\pi bv-\pi bw-\imath b\partial_{u}+\imath b\partial_{w}}).
$$
Substituting all these into the expression for $U\varphi(\mathcal{E}_{1})U^{\ast}$ we obtain
$$
U\varphi(\mathcal{E}_{1})U^{\ast} =
$$
$$
(uv)(vw)e^{-w\partial_{u}}e^{w\partial_{v}}
e^{-\frac{\pi\imath}{4\pi^2 b^{2}}(\pi bu-\pi bv+\pi bw-\imath b\partial_{u}+\imath b\partial_{w})^{2}}
g_{b}(e^{-\pi bu+\pi bv-\pi bw+\imath b\partial_{u}-\imath b\partial_{w}})\times
$$
$$
g_{b}(e^{-2\pi bw})g_{b}(e^{-\pi bu+\pi bv-3\pi bw-\imath b\partial_{u}+\imath b\partial_{w}})g_{b}(e^{-\pi bu+\pi bv-\pi bw-\imath b\partial_{u}+\imath b\partial_{w}})
\varphi(e^{\pi bw-\imath b\partial_{w}})(h.c.).
$$
Note, that $g_{b}(e^{-\pi bu+\pi bv-\pi bw-\imath b\partial_{u}+\imath b\partial_{w}})$ and
$\varphi(e^{\pi bw-\imath b\partial_{w}})$ commute, so the quantum dilogarithm passes through and cancels with its hermitian conjugate:
$$
U\varphi(\mathcal{E}_{1})U^{\ast} =
$$
$$
(uv)(vw)e^{-w\partial_{u}}e^{w\partial_{v}}
e^{-\frac{\pi\imath}{4\pi^2 b^{2}}(\pi bu-\pi bv+\pi bw-\imath b\partial_{u}+\imath b\partial_{w})^{2}}
g_{b}(e^{-\pi bu+\pi bv-\pi bw+\imath b\partial_{u}-\imath b\partial_{w}})\times
$$
$$
g_{b}(e^{-2\pi bw})g_{b}(e^{-\pi bu+\pi bv-3\pi bw-\imath b\partial_{u}+\imath b\partial_{w}})
\varphi(e^{\pi bw-\imath b\partial_{w}})(h.c.).
$$
Let $A$, $B$ be self-adjoint operators satisfying the relation $[A,B] = c$, where $c$ is a number. Let $f(x)$ be a function and $\alpha$ a number. Then
$$
e^{\alpha B^{2}}f(A) = f(A -2c\alpha B)e^{\alpha B^{2}}.
$$
Using this identity to push the exponent $e^{-\frac{\pi\imath}{4\pi^2 b^{2}}(\pi bu-\pi bv+\pi bw-\imath b\partial_{u}+\imath b\partial_{w})^{2}}$ to the right we obtain
$$
U\varphi(\mathcal{E}_{1})U^{\ast} =
(uv)(vw)e^{-w\partial_{u}}e^{w\partial_{v}}
g_{b}(e^{-\pi bu+\pi bv-\pi bw+\imath b\partial_{u}-\imath b\partial_{w}})\times
$$
$$
g_{b}(e^{-\pi bu+\pi bv-3\pi bw+\imath b\partial_{u}-\imath b\partial_{w}})
g_{b}(e^{-2\pi bu+2\pi bv-4\pi bw})
\varphi(e^{\pi bu-\pi bv+2\pi bw-\imath b\partial_{u}})(h.c.)
$$
Now use the relation
$$
e^{\alpha x\partial_{y}}f(y,\partial_{x}) = f(y+\alpha x,\partial_{x}-\alpha\partial_{y})e^{\alpha x\partial_{y}},
$$
to push the exponents $e^{-w\partial_{u}}$ and $e^{w\partial_{v}}$ to the right:
$$
U\varphi(\mathcal{E}_{1})U^{\ast} =
(uv)(vw)g_{b}(e^{-\pi bu+\pi bv+\pi bw+\imath b\partial_{v}-\imath b\partial_{w}})\times
$$
$$
g_{b}(e^{-\pi bu+\pi bv-\pi bw+\imath b\partial_{v}-\imath b\partial_{w}})
g_{b}(e^{-2\pi bu+2\pi bv})
\varphi(e^{\pi bu-\pi bv-\imath b\partial_{u}})(h.c.) =
$$
$$
g_{b}(e^{\pi bu-\pi bv+\pi bw-\imath b\partial_{u}+\imath b\partial_{w}})
g_{b}(e^{-\pi bu-\pi bv+\pi bw-\imath b\partial_{u}+\imath b\partial_{w}})
g_{b}(e^{-2\pi bv+2\pi bw})
\varphi(e^{\pi bv-\pi bw-\imath b\partial_{v}})(h.c.)
$$

Let $\mathcal{E}_{2}$ be the generator in $s_{1}s_{2}s_{1}$ representation. Then
$$
U\varphi(\mathcal{E}_{2})U^{\ast} =
(uv)(vw)e^{-w\partial_{u}}e^{w\partial_{v}}g_{b}(e^{-\pi bu+\pi bv-\pi bw-\imath b\partial_{u}+\imath b\partial_{w}})
g^{\ast}_{b}(e^{\pi bu-\pi bv+\pi bw-\imath b\partial_{u}+\imath b\partial_{w}})\times
$$
$$
g_{b}(e^{\pi bu-\pi bv +\pi bw-\imath b\partial_{u}+\imath b\partial_{w}})
g_{b}(e^{-\pi bu -\pi bv+\pi bw-\imath b\partial_{u}+\imath b\partial_{w}})
g_{b}(e^{-2\pi bv+2\pi bw})
\varphi(e^{\pi bv-\pi bw-\imath b\partial_{v}})\times (h.c.),
$$
where by $(h.c.)$ we denoted the hermitian conjugate operator of everything that stands before  $\varphi(e^{\pi bv-\pi bw-\imath b\partial_{v}})$. Noticing that
$$
g^{\ast}_{b}(e^{\pi bu-\pi bv+\pi bw-\imath b\partial_{u}+\imath b\partial_{w}})
g_{b}(e^{\pi bu-\pi bv +\pi bw-\imath b\partial_{u}+\imath b\partial_{w}}) = 1,
$$
we obtain
$$
U\varphi(\mathcal{E}_{2})U^{\ast} =
(uv)(vw)e^{-w\partial_{u}}e^{w\partial_{v}}\times
$$
$$
g_{b}(e^{-\pi bu+\pi bv-\pi bw-\imath b\partial_{u}+\imath b\partial_{w}})
g_{b}(e^{-\pi bu -\pi bv+\pi bw-\imath b\partial_{u}+\imath b\partial_{w}})
g_{b}(e^{-2\pi bv+2\pi bw})
\varphi(e^{\pi bv-\pi bw-\imath b\partial_{v}})\times (h.c.).
$$
Let $A_{1}$, $A_{2}$ be as follows
$$
A_{1} = e^{-\pi bu+\pi bv-\pi bw-\imath b\partial_{u}+\imath b\partial_{w}},
$$
$$
A_{2} = e^{-2\pi bv+2\pi bw},
$$
Then
$$
q^{-1}A_{1}A_{2} = e^{-\pi bu -\pi bv+\pi bw-\imath b\partial_{u}+\imath b\partial_{w}},
$$
moreover
$$
A_{1}A_{2} = q^{2}A_{2}A_{1},
$$
and we can apply the pentagon identity \cite{Ka0}:
$$
g_{b}(A_{1})g_{b}(q^{-1}A_{1}A_{2})g_{b}(A_{2}) = g_{b}(A_{2})g_{b}(A_{1}),
$$
which leads to the following result
$$
U\varphi(\mathcal{E}_{2})U^{\ast} = (uv)(vw)e^{-w\partial_{u}}e^{w\partial_{v}}
g_{b}(e^{-2\pi bv+2\pi bw})g_{b}(e^{-\pi bu+\pi bv-\pi bw-\imath b\partial_{u}+\imath b\partial_{w}})
\varphi(e^{\pi bv-\pi bw-\imath b\partial_{v}})\times (h.c.).
$$
Operators $g_{b}(e^{-\pi bu+\pi bv-\pi bw-\imath b\partial_{u}+\imath b\partial_{w}})$ and
$\varphi(e^{\pi bv-\pi bw-\imath b\partial_{v}})$ commute, so the quantum dilogarithm passes through and cancels with its conjugate. We obtain
$$
U\varphi(\mathcal{E}_{2})U^{\ast} = (uv)(vw)e^{-w\partial_{u}}e^{w\partial_{v}}
g_{b}(e^{-2\pi bv+2\pi bw})\varphi(e^{\pi bv-\pi bw-\imath b\partial_{v}})\times (h.c.) =
$$
$$
(uv)(vw)g_{b}(e^{-2\pi bv})\varphi(e^{\pi bv-\imath b\partial_{v}})g^{\ast}_{b}(e^{-2\pi bv})(vw)(uv) =
g_{b}(e^{-2\pi bw})\varphi(e^{\pi bw-\imath b\partial_{w}})g^{\ast}_{b}(e^{-2\pi bw}).
$$

The case of $\mathcal{F}_{1}$

$$
U\varphi(\mathcal{F}_{1})U^{\ast} =
$$
$$
(uv)(vw)e^{-w\partial_{u}}e^{w\partial_{v}}g_{b}(e^{-\pi bu+\pi bv-\pi bw-\imath b\partial_{u}+\imath b\partial_{w}})
g^{\ast}_{b}(e^{\pi bu-\pi bv+\pi bw-\imath b\partial_{u}+\imath b\partial_{w}})
g_{b}(e^{\pi bu-\pi bv+\pi bw+\imath b\partial_{u}-\imath b\partial_{w}})\times
$$
$$
g_{b}(e^{-4\pi b\nu_{1}+3\pi bu-\pi bv+\pi bw+\imath b\partial_{u}-\imath b\partial_{w}})
g_{b}(e^{-4\pi b\nu_{1}+4\pi bu-2\pi bv+2\pi bw})
\varphi(e^{2\pi b\nu_{1}-2\pi bu+\pi bv-\pi bw+\imath b\partial_{w}})(h.c.)
$$
The factors $g_{b}(e^{-\pi bu+\pi bv-\pi bw-\imath b\partial_{u}+\imath b\partial_{w}})$ and $g^{\ast}_{b}(e^{\pi bu-\pi bv+\pi bw-\imath b\partial_{u}+\imath b\partial_{w}})$ commute so we can change their order. After that we observe that
$$
g_{b}(e^{-\pi bu+\pi bv-\pi bw-\imath b\partial_{u}+\imath b\partial_{w}})
g_{b}(e^{\pi bu-\pi bv+\pi bw+\imath b\partial_{u}-\imath b\partial_{w}}) =
e^{\frac{\pi\imath}{4\pi^{2}b^{2}}(\pi bu-\pi bv+\pi bw+\imath b\partial_{u}-\imath b\partial_{w})^{2}},
$$
which follows from the identity $g_{b}(x)g_{b}(\frac{1}{x}) = e^{\frac{\pi\imath}{4\pi^{2}b^{2}}\log^{2}x}$. Doing these we obtain
$$
U\varphi(\mathcal{F}_{1})U^{\ast} =
$$
$$
(uv)(vw)e^{-w\partial_{u}}e^{w\partial_{v}}
e^{\frac{\pi\imath}{4\pi^{2}b^{2}}(\pi bu-\pi bv+\pi bw+\imath b\partial_{u}-\imath b\partial_{w})^{2}}
g^{\ast}_{b}(e^{\pi bu-\pi bv+\pi bw-\imath b\partial_{u}+\imath b\partial_{w}})\times
$$
$$
g_{b}(e^{-4\pi b\nu_{1}+3\pi bu-\pi bv+\pi bw+\imath b\partial_{u}-\imath b\partial_{w}})
g_{b}(e^{-4\pi b\nu_{1}+4\pi bu-2\pi bv+2\pi bw})
\varphi(e^{2\pi b\nu_{1}-2\pi bu+\pi bv-\pi bw+\imath b\partial_{w}})(h.c.).
$$
Again using the identity $g_{b}(x)g_{b}(\frac{1}{x}) = e^{\frac{\pi\imath}{4\pi^{2}b^{2}}\log^{2}x}$ to rewrite
$$
g^{\ast}_{b}(e^{\pi bu-\pi bv+\pi bw-\imath b\partial_{u}+\imath b\partial_{w}}) =
e^{-\frac{\pi\imath}{4\pi^{2}b^{2}}(\pi bu-\pi bv+\pi bw-\imath b\partial_{u}+\imath b\partial_{w})^{2}}
g_{b}(e^{-\pi bu+\pi bv-\pi bw+\imath b\partial_{u}-\imath b\partial_{w}}),
$$
we have
$$
U\varphi(\mathcal{F}_{1})U^{\ast} =
$$
$$
(uv)(vw)e^{-w\partial_{u}}e^{w\partial_{v}}
e^{\frac{\pi\imath}{4\pi^{2}b^{2}}(\pi bu-\pi bv+\pi bw+\imath b\partial_{u}-\imath b\partial_{w})^{2}}
e^{-\frac{\pi\imath}{4\pi^{2}b^{2}}(\pi bu-\pi bv+\pi bw-\imath b\partial_{u}+\imath b\partial_{w})^{2}}
g_{b}(e^{-\pi bu+\pi bv-\pi bw+\imath b\partial_{u}-\imath b\partial_{w}})\times
$$
$$
g_{b}(e^{-4\pi b\nu_{1}+3\pi bu-\pi bv+\pi bw+\imath b\partial_{u}-\imath b\partial_{w}})
g_{b}(e^{-4\pi b\nu_{1}+4\pi bu-2\pi bv+2\pi bw})
\varphi(e^{2\pi b\nu_{1}-2\pi bu+\pi bv-\pi bw+\imath b\partial_{w}})(h.c.).
$$
Let
$$
A_{1} = e^{-\pi bu+\pi bv-\pi bw+\imath b\partial_{u}-\imath b\partial_{w}},
$$
$$
A_{2} = e^{-4\pi b\nu_{1}+4\pi bu-2\pi bv+2\pi bw}.
$$
Then
$$
q^{-1}A_{1}A_{2} = e^{-4\pi b\nu_{1}+3\pi bu-\pi bv+\pi bw+\imath b\partial_{u}-\imath b\partial_{w}},
$$
$$
A_{1}A_{2} = q^{2}A_{2}A_{1},
$$
and we can apply the pentagon identity $g_{b}(A_{1})g_{b}(q^{-1}A_{1}A_{2})g_{b}(A_{2}) = g_{b}(A_{2})g_{b}(A_{1})$:
$$
g_{b}(e^{-\pi bu+\pi bv-\pi bw+\imath b\partial_{u}-\imath b\partial_{w}})
g_{b}(e^{-4\pi b\nu_{1}+3\pi bu-\pi bv+\pi bw+\imath b\partial_{u}-\imath b\partial_{w}})
g_{b}(e^{-4\pi b\nu_{1}+4\pi bu-2\pi bv+2\pi bw}) =
$$
$$
g_{b}(e^{-4\pi b\nu_{1}+4\pi bu-2\pi bv+2\pi bw})g_{b}(e^{-\pi bu+\pi bv-\pi bw+\imath b\partial_{u}-\imath b\partial_{w}})
$$
Note also that $g_{b}(e^{-\pi bu+\pi bv-\pi bw+\imath b\partial_{u}-\imath b\partial_{w}})$ commutes with $\varphi(e^{2\pi b\nu_{1}-2\pi bu+\pi bv-\pi bw+\imath b\partial_{w}})$ so we obtain
$$
U\varphi(\mathcal{F}_{1})U^{\ast} =
$$
$$
(uv)(vw)e^{-w\partial_{u}}e^{w\partial_{v}}
e^{\frac{\pi\imath}{4\pi^{2}b^{2}}(\pi bu-\pi bv+\pi bw+\imath b\partial_{u}-\imath b\partial_{w})^{2}}
e^{-\frac{\pi\imath}{4\pi^{2}b^{2}}(\pi bu-\pi bv+\pi bw-\imath b\partial_{u}+\imath b\partial_{w})^{2}}\times
$$
$$
g_{b}(e^{-4\pi b\nu_{1}+4\pi bu-2\pi bv+2\pi bw})
\varphi(e^{2\pi b\nu_{1}-2\pi bu+\pi bv-\pi bw+\imath b\partial_{w}})(h.c.)
$$
To push the quadratic exponents to the right we use the formula $e^{\alpha B^{2}}f(A) = f(A -2c\alpha B)e^{\alpha B^{2}}$ for self-adjoint $A$ and $B$ satisfying the relation $[A,B] = c$, where $c$, $\alpha$ are numbers:
$$
U\varphi(\mathcal{F}_{1})U^{\ast} =
(uv)(vw)e^{-w\partial_{u}}e^{w\partial_{v}}
e^{\frac{\pi\imath}{4\pi^{2}b^{2}}(\pi bu-\pi bv+\pi bw+\imath b\partial_{u}-\imath b\partial_{w})^{2}}\times
$$
$$
g_{b}(e^{-4\pi b\nu_{1} +3\pi bu-\pi bv+\pi bw+\imath b\partial_{u}-\imath b\partial_{w}})
\varphi(e^{2\pi b\nu_{1}-2\pi bu+\pi bv-\pi bw+\imath b\partial_{w}})(h.c.) =
$$
$$
(uv)(vw)e^{-w\partial_{u}}e^{w\partial_{v}}
g_{b}(e^{-4\pi b\nu_{1}+2\pi bu})
\varphi(e^{2\pi b\nu_{1}-\pi bu+\imath b\partial_{u}})(h.c.) =
$$
$$
(uv)(vw)
g_{b}(e^{-4\pi b\nu_{1}+2\pi bu-2\pi bw})
\varphi(e^{2\pi b\nu_{1}-\pi bu+\pi bw+\imath b\partial_{u}})(h.c.) =
$$
$$
g_{b}(e^{-4\pi b\nu_{1}-2\pi bu+2\pi bv})
\varphi(e^{2\pi b\nu_{1}+\pi bu-\pi bv+\imath b\partial_{v}})(h.c.).
$$

$\Box$

Let $K_{1}$, $\mathcal{E}_{1}$, $\mathcal{F}_{1}$ be the subset of $U_{q}(\mathfrak{sl}(3))$ generators in $s_{1}s_{2}s_{1}$ principal series representation. Recall that for a complex-valued function $\varphi(x)$ we have

$$
\varphi(K_{1}) = \varphi(e^{-2\pi b\nu_{1}+2\pi bu-\pi bv+2\pi bw}),
$$

$$
\varphi(\mathcal{E}_{1}) = g_{b}(e^{-2\pi bw})\varphi(e^{\pi bw-\imath b\partial_{w}})g_{b}^{\ast}(e^{-2\pi bw}),
$$

$$
\varphi(\mathcal{F}_{1}) = g_{b}(e^{\pi bu-\pi bv+\pi bw+\imath b\partial_{u}-\imath b\partial_{w}})
g_{b}(e^{-4\pi b\nu_{1}+3\pi bu-\pi bv+\pi bw+\imath b\partial_{u}-\imath b\partial_{w}})
g_{b}(e^{-4\pi b\nu_{1}+4\pi bu-2\pi bv+2\pi bw})\times
$$
$$
\varphi(e^{2\pi b\nu_{1}-2\pi bu+\pi bv-\pi bw+\imath b\partial_{w}})\times
$$
$$
g_{b}^{\ast}(e^{-4\pi b\nu_{1}+4\pi bu-2\pi bv+2\pi bw})
g_{b}^{\ast}(e^{-4\pi b\nu_{1}+3\pi bu-\pi bv+\pi bw+\imath b\partial_{u}-\imath b\partial_{w}})
g_{b}^{\ast}(e^{\pi bu-\pi bv+\pi bw+\imath b\partial_{u}-\imath b\partial_{w}}),
$$

This subset generates $U_{q}(\mathfrak{sl}(2))$ subalgebra.

The unitary transform from the following proposition was given in the proof of Theorem 4.7 in \cite{Ip1}. It was used as the first step in mapping the generators $K_{i}$, $\mathcal{E}_{i}$, $\mathcal{F}_{i}$ of $U_{q}(\mathfrak{sl}(2))_{i}$ subalgebra of $U_{q}(\mathfrak{g})$ to the formulas corresponding to positive principal series representations of $U_{q}(\mathfrak{sl}(2))$.  We explicitly check its action on the functions of generators.

\begin{lem}
Let $K_{1}$, $\mathcal{E}_{1}$, $\mathcal{F}_{1}$ be as above. Let $\varphi(x)$ be a complex-valued function. Let
\begin{equation}
V = e^{(\nu_{1}+\frac{v}{2})\partial_{w}}e^{(\nu_{1}+\frac{v}{2})\partial_{u}}
e^{-\frac{\pi\imath u^{2}}{2}+2\pi\imath\nu_{1}u}e^{-u\partial_{w}}e^{-\frac{\pi\imath w^{2}}{2}}g_{b}(e^{2\pi bw})
g_{b}(e^{4\pi b\nu_{1}-2\pi bu})
\end{equation}
be a unitary transform. Then
\begin{equation}
V\varphi(K_{1})V^{\ast} = \varphi(e^{2\pi bw}),
\end{equation}

\begin{equation}
V\varphi(\mathcal{E}_{1})V^{\ast} = \varphi(e^{-\imath b\partial_{w}}),
\end{equation}

\begin{multline}
V\varphi(\mathcal{F}_{1})V^{\ast} = g_{b}(e^{2\pi bw}e^{-2\pi bu})g_{b}(e^{2\pi bw}e^{\imath b\partial_{u}})g_{b}(e^{2\pi bw}e^{2\pi bu})\times     \\
\varphi(e^{-2\pi bw+\imath b\partial_{w}})\times  \\
g^{\ast}_{b}(e^{2\pi bw}e^{2\pi bu})g^{\ast}_{b}(e^{2\pi bw}e^{\imath b\partial_{u}})g^{\ast}_{b}(e^{2\pi bw}e^{-2\pi bu}).
\end{multline}
\end{lem}
$\proof$
$$
V\varphi(\mathcal{F}_{1})V^{\ast} =
e^{(\nu_{1}+\frac{v}{2})\partial_{w}}e^{(\nu_{1}+\frac{v}{2})\partial_{u}}
e^{-\frac{\pi\imath u^{2}}{2}+2\pi\imath\nu_{1}u}e^{-u\partial_{w}}e^{-\frac{\pi\imath w^{2}}{2}}g_{b}(e^{2\pi bw})\times
$$
$$
g_{b}(e^{4\pi b\nu_{1}-2\pi bu})
g_{b}(e^{\pi bu-\pi bv+\pi bw+\imath b\partial_{u}-\imath b\partial_{w}})
g_{b}(e^{-4\pi b\nu_{1}+3\pi bu-\pi bv+\pi bw+\imath b\partial_{u}-\imath b\partial_{w}})
g_{b}(e^{-4\pi b\nu_{1}+4\pi bu-2\pi bv+2\pi bw})\times
$$
$$
\varphi(e^{2\pi b\nu_{1}-2\pi bu+\pi bv-\pi bw+\imath b\partial_{w}})\times (h.c.).
$$
Let
$$
U_{1} = e^{4\pi b\nu_{1}-2\pi bu},
$$
$$
U_{2} = e^{-4\pi b\nu_{1}+3\pi bu-\pi bv+\pi bw+\imath b\partial_{u}-\imath b\partial_{w}}.
$$
Then
$$
q^{-1}U_{1}U_{2} = e^{\pi bu-\pi bv+\pi bw+\imath b\partial_{u}-\imath b\partial_{w}},
$$
and the following relation holds
$$
U_{1}U_{2} = q^{2}U_{2}U_{1},
$$
which allows us to use the quantum pentagon identity
$$
g_{b}(U_{1})g_{b}(q^{-1}U_{1}U_{2})g_{b}(U_{2}) = g_{b}(U_{2})g_{b}(U_{1}).
$$
We obtain
$$
V\varphi(\mathcal{F}_{1})V^{\ast} =
e^{(\nu_{1}+\frac{v}{2})\partial_{w}}e^{(\nu_{1}+\frac{v}{2})\partial_{u}}
e^{-\frac{\pi\imath u^{2}}{2}+2\pi\imath\nu_{1}u}e^{-u\partial_{w}}e^{-\frac{\pi\imath w^{2}}{2}}g_{b}(e^{2\pi bw})\times
$$
$$
g_{b}(e^{-4\pi b\nu_{1}+3\pi bu-\pi bv+\pi bw+\imath b\partial_{u}-\imath b\partial_{w}})
g_{b}(e^{4\pi b\nu_{1}-2\pi bu})
g_{b}(e^{-4\pi b\nu_{1}+4\pi bu-2\pi bv+2\pi bw})
\varphi(e^{2\pi b\nu_{1}-2\pi bu+\pi bv-\pi bw+\imath b\partial_{w}})\times (h.c.)
$$
Note that the operator $g_{b}(e^{4\pi b\nu_{1}-2\pi bu})$ commutes with $g_{b}(e^{-4\pi b\nu_{1}+4\pi bu-2\pi bv+2\pi bw})$ and $\varphi(e^{2\pi b\nu_{1}-2\pi bu+\pi bv-\pi bw+\imath b\partial_{w}})$, so it goes through and cancels with its hermitian conjugate.
$$
V\varphi(\mathcal{F}_{1})V^{\ast} =
e^{(\nu_{1}+\frac{v}{2})\partial_{w}}e^{(\nu_{1}+\frac{v}{2})\partial_{u}}
e^{-\frac{\pi\imath u^{2}}{2}+2\pi\imath\nu_{1}u}e^{-u\partial_{w}}e^{-\frac{\pi\imath w^{2}}{2}}g_{b}(e^{2\pi bw})\times
$$
$$
g_{b}(e^{-4\pi b\nu_{1}+3\pi bu-\pi bv+\pi bw+\imath b\partial_{u}-\imath b\partial_{w}})
g_{b}(e^{-4\pi b\nu_{1}+4\pi bu-2\pi bv+2\pi bw})
\varphi(e^{2\pi b\nu_{1}-2\pi bu+\pi bv-\pi bw+\imath b\partial_{w}})\times (h.c.).
$$
Using the following operator relations
$$
e^{\alpha x^{2} +\beta x}f(\partial_{x}) = f(\partial_{x}-2\alpha x -\beta)e^{\alpha x^{2} +\beta x},
$$
$$
e^{\alpha x\partial_{y}}f(y,\partial_{x}) = f(y+\alpha x,\partial_{x}-\alpha\partial_{y})e^{\alpha x\partial_{y}},
$$
we move all the exponents to the right
$$
V\varphi(\mathcal{F}_{1})V^{\ast} =
e^{(\nu_{1}+\frac{v}{2})\partial_{w}}e^{(\nu_{1}+\frac{v}{2})\partial_{u}}
e^{-\frac{\pi\imath u^{2}}{2}+2\pi\imath\nu_{1}u}e^{-u\partial_{w}}g_{b}(e^{2\pi bw})\times
$$
$$
g_{b}(e^{-4\pi b\nu_{1}+3\pi bu-\pi bv+2\pi bw+\imath b\partial_{u}-\imath b\partial_{w}})
g_{b}(e^{-4\pi b\nu_{1}+4\pi bu-2\pi bv+2\pi bw})
\varphi(e^{2\pi b\nu_{1}-2\pi bu+\pi bv-2\pi bw+\imath b\partial_{w}})\times (h.c.) =
$$
$$
e^{(\nu_{1}+\frac{v}{2})\partial_{w}}e^{(\nu_{1}+\frac{v}{2})\partial_{u}}
e^{-\frac{\pi\imath u^{2}}{2}+2\pi\imath\nu_{1}u}g_{b}(e^{2\pi bw-2\pi bu})\times
$$
$$
g_{b}(e^{-4\pi b\nu_{1}+\pi bu-\pi bv+2\pi bw+\imath b\partial_{u}})
g_{b}(e^{-4\pi b\nu_{1}+2\pi bu-2\pi bv+2\pi bw})
\varphi(e^{2\pi b\nu_{1}+\pi bv-2\pi bw+\imath b\partial_{w}})\times (h.c.) =
$$
$$
e^{(\nu_{1}+\frac{v}{2})\partial_{w}}e^{(\nu_{1}+\frac{v}{2})\partial_{u}}
g_{b}(e^{2\pi bw-2\pi bu})
g_{b}(e^{-2\pi b\nu_{1}-\pi bv+2\pi bw+\imath b\partial_{u}})\times
$$
$$
g_{b}(e^{-4\pi b\nu_{1}+2\pi bu-2\pi bv+2\pi bw})
\varphi(e^{2\pi b\nu_{1}+\pi bv-2\pi bw+\imath b\partial_{w}})\times (h.c.) =
$$
$$
g_{b}(e^{2\pi bw-2\pi bu})g_{b}(e^{2\pi bw+\imath b\partial_{u}})
g_{b}(e^{2\pi bw+2\pi bu})
\varphi(e^{-2\pi bw+\imath b\partial_{w}})
g^{\ast}_{b}(e^{2\pi bw+2\pi bu})g^{\ast}_{b}(e^{2\pi bw+\imath b\partial_{u}})
g^{\ast}_{b}(e^{2\pi bw-2\pi bu}).
$$

$\Box$

To finish the mapping $\varphi(K_{1})$,$\varphi(\mathcal{E}_{1})$,$\varphi(\mathcal{F}_{1})\rightarrow$ $\varphi(K)$,$\varphi(\mathcal{E})$,$\varphi(\mathcal{F})$,  where the second set of operators is defined by the equations (\ref{function of E})-(\ref{function of F}), we need to perform a certain integral transformation which will be defined shortly.

Let $\lambda$ be a positive real number. Define the following set of functions, \cite{Ka2}, \cite{KLSTS}
\begin{equation}
\Phi_{\lambda}(u) = e^{\pi\imath u^{2}+\pi Qu}G_{b}(-\imath u+\imath\lambda)G_{b}(-\imath u-\imath\lambda).
\end{equation}
The integral transform $\Phi$ is defined by
$$
\Phi: L^{2}(\mathbb{R}) \rightarrow L^{2}(\mathbb{R}^{+},d\mu(\lambda)),
$$
\begin{equation}\label{Kashaev integral transform}
\Phi: f(u) \rightarrow F(\lambda) = \int\limits_{\mathbb{R}-\imath 0} du f(u)\Phi^{\ast}_{\lambda}(u),
\end{equation}
This transform is an isometry, see \cite{Ka2}. The inverse is given by
$$
\Phi^{-1}: L^{2}(\mathbb{R}^{+},d\mu(\lambda)) \rightarrow L^{2}(\mathbb{R}),
$$
\begin{equation}\label{Kashaev inverse integral transform}
\Phi^{-1} : F(\lambda) \rightarrow f(u) = \lim\limits_{\epsilon\rightarrow 0}\int\limits_{0}^{+\infty}F(\lambda)\Phi_{\lambda}(u+\imath\epsilon)e^{-2\pi\epsilon u} d\mu(\lambda),
\end{equation}
with the measure given by $d\mu(\lambda) = 4\sinh(\pi b\lambda)\sinh(\pi b^{-1}\lambda)$.

\begin{prop}
The function $\Phi_{\lambda}(u)$ is an eigenfunction of the operator $g_{b}(e^{2\pi bw}e^{-2\pi bu})g_{b}(e^{2\pi bw}e^{\imath b\partial_{u}})g_{b}(e^{2\pi bw}e^{2\pi bu})$:
\begin{equation}
g_{b}(e^{2\pi bw}e^{-2\pi bu})g_{b}(e^{2\pi bw}e^{\imath b\partial_{u}})g_{b}(e^{2\pi bw}e^{2\pi bu})\Phi_{\lambda}(u) =
g_{b}(e^{2\pi b\lambda+2\pi bw})g_{b}(e^{-2\pi b\lambda+2\pi bw})\Phi_{\lambda}(u).
\end{equation}
\end{prop}
$\proof$

Recalling the definition of $g_{b}(x)$
$$
g_{b}(x) = \frac{\bar{\zeta}_{b}}{G_{b}(\frac{Q}{2} + \frac{1}{2\pi\imath b}\log x)},
$$
and the Fourier transform
$$
g_{b}(x) = \int d\tau x^{\imath b^{-1}\tau}e^{\pi Q\tau}G_{b}(-\imath\tau),
$$
we obtain
$$
g_{b}(e^{2\pi bw}e^{-2\pi bu}) = \frac{\bar{\zeta}_{b}}{G_{b}(\frac{Q}{2} -\imath w + \imath u)},
$$
$$
g_{b}(e^{2\pi bw}e^{2\pi bu}) = \frac{\bar{\zeta}_{b}}{G_{b}(\frac{Q}{2} -\imath w - \imath u)},
$$
$$
g_{b}(e^{2\pi bw}e^{\imath b\partial_{u}}) = \int d\tau e^{\pi Q\tau+2\pi\imath w\tau}G_{b}(-\imath\tau)e^{-\tau\partial_{u}}.
$$
Substituting these expressions into the left-hand side of the eigenvalue equation we obtain
$$
g_{b}(e^{2\pi bw}e^{-2\pi bu})g_{b}(e^{2\pi bw}e^{\imath b\partial_{u}})g_{b}(e^{2\pi bw}e^{2\pi bu})\Phi_{\lambda}(u) =
$$
$$
\frac{1}{G_{b}(\frac{Q}{2}-\imath w+\imath u)}\int d\tau e^{\pi Q\tau + 2\pi\imath w\tau}G_{b}(-\imath\tau)e^{-\tau\partial_{u}}
\frac{e^{\pi\imath u^{2}+\pi Qu}G_{b}(-\imath u+\imath\lambda)G_{b}(-\imath u-\imath\lambda)}{G_{b}(\frac{Q}{2}-\imath w-\imath u)} =
$$
$$
\frac{1}{G_{b}(\frac{Q}{2}-\imath w+\imath u)}\int d\tau e^{\pi Q\tau+2\pi\imath w\tau+\pi\imath(u-\tau)^{2}+\pi Q(u-\tau)}
\frac{G_{b}(-\imath\tau)G_{b}(-\imath u+\imath\lambda+\imath\tau)G_{b}(-\imath u-\imath\lambda+\imath\tau)}{G_{b}(\frac{Q}{2}-\imath w-\imath u+\imath\tau)}
$$
Applying the reflection formula
$$
G_{b}(-\imath\tau) = \frac{e^{-\pi\imath\tau^{2}-\pi Q\tau}}{G_{b}(Q+\imath\tau)},
$$
we get
$$
g_{b}(e^{2\pi bw}e^{-2\pi bu})g_{b}(e^{2\pi bw}e^{\imath b\partial_{u}})g_{b}(e^{2\pi bw}e^{2\pi bu})\Phi_{\lambda}(u) =
$$
$$
\frac{e^{\pi\imath u^{2}+\pi Qu}}{G_{b}(\frac{Q}{2}-\imath w+\imath u)}\int d\tau
e^{-2\pi(\frac{Q}{2}+\imath u-\imath w)\tau}
\frac{G_{b}(-\imath u+\imath\lambda+\imath\tau)G_{b}(-\imath u-\imath\lambda+\imath\tau)}{G_{b}(\frac{Q}{2}-\imath w-\imath u+\imath\tau)G_{b}(Q+\imath\tau)}.
$$
Let $\alpha = -\imath u+\imath\lambda$, $\beta = -\imath u-\imath\lambda$, $\gamma = \frac{Q}{2}+\imath u-\imath w$. Then
$$
g_{b}(e^{2\pi bw}e^{-2\pi bu})g_{b}(e^{2\pi bw}e^{\imath b\partial_{u}})g_{b}(e^{2\pi bw}e^{2\pi bu})\Phi_{\lambda}(u) =
$$
$$
\frac{e^{\pi\imath u^{2}+\pi Qu}}{G_{b}(\frac{Q}{2}-\imath w+\imath u)}\int d\tau
e^{-2\pi\gamma\tau}
\frac{G_{b}(\alpha+\imath\tau)G_{b}(\beta+\imath\tau)}{G_{b}(\alpha+\beta+\gamma+\imath\tau)G_{b}(Q+\imath\tau)} =
$$
$$
\frac{e^{\pi\imath u^{2}+\pi Qu}}{G_{b}(\frac{Q}{2}-\imath w+\imath u)}
\frac{G_{b}(\alpha)G_{b}(\beta)G_{b}(\gamma)}{G_{b}(\alpha+\gamma)G_{b}(\beta+\gamma)},
$$
here $4-5$ integral identity \cite{V} has been used. Substituting $\alpha$, $\beta$, $\gamma$ we obtain
$$
g_{b}(e^{2\pi bw}e^{-2\pi bu})g_{b}(e^{2\pi bw}e^{\imath b\partial_{u}})g_{b}(e^{2\pi bw}e^{2\pi bu})\Phi_{\lambda}(u) =
e^{\pi\imath u^{2}+\pi Qu}\frac{G_{b}(-\imath u+\imath\lambda)G_{b}(-\imath u-\imath\lambda)}{G_{b}(\frac{Q}{2}+\imath\lambda-\imath w)G_{b}(\frac{Q}{2}-\imath\lambda-\imath w)} =
$$
$$
\frac{1}{G_{b}(\frac{Q}{2}+\imath\lambda-\imath w)G_{b}(\frac{Q}{2}-\imath\lambda-\imath w)}\Phi_{\lambda}(u) =
g_{b}(e^{2\pi b\lambda+2\pi bw})g_{b}(e^{-2\pi b\lambda+2\pi bw})\Phi_{\lambda}(u).
$$
$\Box$

\begin{cor}
Let $\varphi(K_{1})$,$\varphi(\mathcal{E}_{1})$,$\varphi(\mathcal{F}_{1})$ be the functions of the subset of generators of $U_{q}(\mathfrak{sl}(3))$ in the positive principal series representation corresponding to the reduced expression $w_{0} = s_{1}s_{2}s_{1}$ of the longest Weyl element. Let $\Phi$, $\Phi^{\ast}$ be the integral transform and its inverse defined in (\ref{Kashaev integral transform})-(\ref{Kashaev inverse integral transform}). Let $\Omega$ be the unitary transform defined by
\begin{equation}
    \Omega_{1} =
    e^{-\frac{\pi\imath}{2}(\lambda+w)^{2}}g_{b}(e^{-2\pi b\lambda-2\pi bw})
    \circ\Phi \circ e^{(\nu_{1}+\frac{v}{2})\partial_{w}}e^{(\nu_{1}+\frac{v}{2})\partial_{u}}
e^{-\frac{\pi\imath u^{2}}{2}+2\pi\imath\nu_{1}u}e^{-u\partial_{w}}e^{-\frac{\pi\imath w^{2}}{2}}g_{b}(e^{2\pi bw})
g_{b}(e^{4\pi b\nu_{1}-2\pi bu}).
\end{equation}
Then
\begin{equation}
\Omega_{1}\varphi(K_{1})\Omega_{1}^{\ast} = \varphi(e^{2\pi bw}),
\end{equation}
\begin{equation}
\Omega_{1}\varphi(\mathcal{E}_{1})\Omega_{1}^{\ast} = g_{b}(e^{-2\pi b\lambda-2\pi bw})\varphi(e^{\pi b\lambda +\pi bw - \imath b\partial_{w}})g^{\ast}_{b}(e^{-2\pi b\lambda-2\pi bw}),
\end{equation}
\begin{equation}
\Omega_{1}\varphi(\mathcal{F}_{1})\Omega_{1}^{\ast} = g_{b}(e^{-2\pi b\lambda+2\pi bw})\varphi(e^{\pi b\lambda-\pi bw+\imath b\partial_{w}})g^{\ast}_{b}(e^{-2\pi b\lambda+2\pi bw}).
\end{equation}
\end{cor}

\begin{cor}
Let $\varphi(K_{2})$,$\varphi(\mathcal{E}_{2})$,$\varphi(\mathcal{F}_{2})$ be the functions of the subset of generators of $U_{q}(\mathfrak{sl}(3))$ in the positive principal series representation corresponding to the reduced expression $w_{0} = s_{2}s_{1}s_{2}$ of the longest Weyl element. Let $\Omega_{2}$ be a unitary transform defined by
\begin{equation}
    \Omega_{2} =
    e^{-\frac{\pi\imath}{2}(\lambda+w)^{2}}g_{b}(e^{-2\pi b\lambda-2\pi bw})
    \circ\Phi \circ e^{(\nu_{2}+\frac{v}{2})\partial_{w}}e^{(\nu_{2}+\frac{v}{2})\partial_{u}}
e^{-\frac{\pi\imath u^{2}}{2}+2\pi\imath\nu_{1}u}e^{-u\partial_{w}}e^{-\frac{\pi\imath w^{2}}{2}}g_{b}(e^{2\pi bw})
g_{b}(e^{4\pi b\nu_{2}-2\pi bu}).
\end{equation}
Then
\begin{equation}
\Omega_{2}\varphi(K_{2})\Omega_{2}^{\ast} = \varphi(e^{2\pi bw}),
\end{equation}
\begin{equation}
\Omega_{2}\varphi(\mathcal{E}_{2})\Omega_{2}^{\ast} = g_{b}(e^{-2\pi b\lambda-2\pi bw})\varphi(e^{\pi b\lambda +\pi bw - \imath b\partial_{w}})g^{\ast}_{b}(e^{-2\pi b\lambda-2\pi bw}),
\end{equation}
\begin{equation}
\Omega_{2}\varphi(\mathcal{F}_{2})\Omega_{2}^{\ast} = g_{b}(e^{-2\pi b\lambda+2\pi bw})\varphi(e^{\pi b\lambda-\pi bw+\imath b\partial_{w}})g^{\ast}_{b}(e^{-2\pi b\lambda+2\pi bw}).
\end{equation}
\end{cor}
\proof
Note, that swapping the indices $K_{1}\leftrightarrow K_{2}$, $\mathcal{E}_{1}\leftrightarrow \mathcal{E}_{2}$, $\mathcal{F}_{1}\leftrightarrow \mathcal{F}_{2}$, $\nu_{1}\leftrightarrow \nu_{2}$, of generators and parameters in a representation of $U_{q}(\mathfrak{sl}(3))$ corresponding to a particular choice of reduced expression of the longest Weyl element gives representation for another choice of reduced expression. So, from the statement, that $\Omega_{1}$ transforms the action of operators
$\varphi(K_{1})$,$\varphi(\mathcal{E}_{1})$,$\varphi(\mathcal{F}_{1})$ in $s_{1}s_{2}s_{1}$ representation to the $U_{q}(\mathfrak{sl}(2))$ formulas (\ref{function of E})-(\ref{function of F}), it automatically follows that $\Omega_{2}$ which is obtained from $\Omega_{1}$ by the replacement of the parameter $\nu_{1}$ by $\nu_{2}$, transforms the action of operators
$\varphi(K_{2})$,$\varphi(\mathcal{E}_{2})$,$\varphi(\mathcal{F}_{2})$ in $s_{2}s_{1}s_{2}$ representation to  $U_{q}(\mathfrak{sl}(2))$ formulas.
$\Box$

\begin{cor}
In the positive principal series representation corresponding to any reduced expression of the Weyl element the generalized Kac's identity holds.
\end{cor}
\proof
As follows from the corollaries 4.2, 4.3, there is a unitary transformation which transforms the operators $\varphi(K_{i})$,$\varphi(\mathcal{E}_{i})$,$\varphi(\mathcal{F}_{i})$ defined in the positive principal series of $U_{q}(\mathfrak{sl}(3))$ to the operators $\varphi(K)$,$\varphi(\mathcal{E})$,$\varphi(\mathcal{F})$ defined in positive principal series representation of $U_{q}(\mathfrak{sl}(2))$. Since the generalized Kac's identity is valid in $U_{q}(\mathfrak{sl}(2))$ case, it follows that it is as well valid in the case of $U_{q}(\mathfrak{sl}(3))$.
$\Box$\\*
This completes the proof of the Theorem 4.1.

Note, that we have also given another proof in $U_{q}(\mathfrak{sl}(3))$ for the Theorem 4.7 in \cite{Ip1} which states that
the positive principal series representation of $U_{q}(\mathfrak{sl}(3))$ decomposes into direct integral of positive principal series representations of its $U_{q}(\mathfrak{sl}(2))$ subalgebra corresponding to any simple root.

\end{document}